\documentclass[aos,seceqn,citesort,dvips]{arximspdf}

\usepackage{graphics}
\doi{10.1214/07-AOS563}
\volume{37}
\issue{1}
\pubyear{2009}
\firstpage{35}
\lastpage{72}

\makeatletter

\newtheorem{theorem}{Theorem}
\newtheorem{lemma}{Lemma}
\newtheorem{proposition}{Proposition}
\newproclaim{remark}{Remark}

\makeatother

\begin{document}
\begin{frontmatter}

\title{Smoothing splines estimators for functional linear~regression}
\runtitle{Functional linear regression}

\begin{aug}
\author[A]{\fnms{Christophe} \snm{Crambes}\ead[label=e1]{Christophe.Crambes@math.ups-tlse.fr}},
\author[B]{\fnms{Alois} \snm{Kneip}\corref{}\ead[label=e2]{akneip@uni-bonn.de}}
and
\author[A]{\fnms{Pascal} \snm{Sarda}\ead[label=e3]{Pascal.Sarda@math.ups-tlse.fr}}
\runauthor{C. Crambes, A. Kneip and P. Sarda}
\affiliation{Universit\'e Paul Sabatier, Universit\'e Paul Sabatier and
Universit\"at Bonn}
\address[A]{C. Crambes\\
P. Sarda\\
Universit\'e Paul Sabatier\\
Institut de Math\'ematiques\\
UMR 5219\\
Laboratoire de Statistique et Probabilit\'es\\
118 Route de Narbonne\\
31062 Toulouse Cedex\\
France\\
\printead{e1}\\
\phantom{E-mail: }\printead*{e3}}
\address[B]{A. Kneip\\
Statistische Abteilung\\
Department of Economics and Hausdorff\\
Center for Mathematics\\
Universit\"{a}t Bonn\\
Adenauerallee 24-26\\
53113 Bonn\\
Germany\\
\printead{e2}}
\end{aug}

\received{\smonth{5} \syear{2007}}
\revised{\smonth{10} \syear{2007}}

%
\begin{abstract}
The paper considers functional linear regression, where scalar responses
$Y_1,\ldots,Y_n$ are modeled in dependence of
random functions $X_1,\ldots,X_n$.
We propose a smoothing splines estimator for the functional slope parameter
based on a
slight modification of the usual penalty. Theoretical analysis
concentrates on
the error in an out-of-sample prediction of the response for a new
random function $X_{n+1}$.
It is shown that rates of convergence of the prediction error depend on
the smoothness of the slope function and on the structure of the predictors.
We then prove that these rates are optimal in the sense that they are minimax
over large classes of possible slope functions and distributions of the
predictive curves. For the case of models with errors-in-variables the
smoothing spline estimator
is modified by
using a denoising correction of the covariance matrix of discretized curves.
The methodology is then applied to a real case study where the aim is
to predict the maximum of the
concentration of
ozone by using the curve of this concentration measured the preceding day.
\end{abstract}

%
\pdfkeywords{Functional linear regression,
functional parameter,
functional variable,
smoothing splines}

\end{frontmatter}

\section{Introduction} \label{sec1}

In a number of important applications the outcome of a response
variable $Y$ depends on the variation of
an explanatory variable $X$ over time (or age, etc.). An example is the
application motivating our study: the data
consist in repeated measurements of pollutant indicators in the area of Toulouse
over the course of a day that are used to explain
the maximum (peak) of pollution for the next day. Generally, a linear
regression model linking observations
$Y_i$ of a response variable with $p$ repeated measures of an
explanatory variable may be written in the form
%
\begin{equation}\label{funclinearregressiondiscr}
Y_{i} = \alpha_0+\frac{1}{p} \sum_{j=1}^{p} \alpha_j
X_{i} (t_{j}) + \varepsilon_{i}^*,\qquad  i=1,\ldots,n.
\end{equation}
Here $t_{1} < \cdots < t_{p}$ denote observation points which are
assumed to belong to a compact interval $I\subset\mathbb{R}$.
The possibly varying strength of the influence of $X_i$ at each
measurement point $t_j$ is quantified by different coefficients
$\alpha_j$. Frequently $p\gg n$ and/or there is a high degree of
collinearity between the ``predictors'' $X_i(t_j),j=1,\ldots,p$, and
standard regression methods are not applicable. In addition,
(\ref{funclinearregressiondiscr}) may incorporate a discretization
error, since one will often have to assume that $Y_i$
also depends on unobserved time points $t$ in between the observation
times $t_j$.
As pointed out by several
authors (Marx and Eilers \cite{MaEi99}, Ramsay and Silverman \cite
{RaSi05} or Cuevas, Febrero and Fraiman \cite{CuFeFr02}) the use of
functional models for these settings has some advantages over discrete,
multivariate approaches. Only in
a functional framework is it possible to profit from qualitative
assumptions like smoothness of
underlying curves. Assuming square integrable functions $X_i$ on
$I\subset\mathbb{R}$, the basic object of our study is a \textit
{functional linear regression model}
%
\begin{equation}\label{funclinearregression}
Y_i = \alpha_0+\int_{I} \alpha(t) X_i(t) \,dt + \varepsilon_i,\qquad
i=1,\ldots,n,
\end{equation}
where $\varepsilon_i$'s are i.i.d. centered random errors,
$\mathbb{E} (\varepsilon_i) = 0$, with variance
$\mathbb{E} (\varepsilon_i^{2}) = \sigma_{\varepsilon}^{2}$, and $\alpha$ is a
square integrable functional parameter defined on $I$ that must be
estimated from the pairs $(X_i,Y_i), i=1,\ldots,n$. This type of
regression model was first considered in Ramsay and Dalzell \cite{RaDa91}.
Obviously, (\ref{funclinearregression}) constitutes
a continuous version of (\ref{funclinearregressiondiscr}), and both
models are linked by
%
\begin{equation}
\varepsilon_i^*=d_i+\varepsilon_i,\qquad \mbox{where }
d_{i} = \int_{I} \alpha(t) X_{i} (t) \,dt - \frac{1}{p} \sum_{j=1}^{p}
\alpha(t_{j}) X_{i} (t_{j})
\end{equation}
may be interpreted as a discretization error, and $\alpha(t_j)=\alpha_j$.

As a consequence of developments of modern technology, data that may
be described by functional regression models can be found in a
lot of fields such as medicine, linguistics, chemometrics
(see, e.g., Ramsay and Silverman \cite{RaSi02,RaSi05} and
Ferraty and Vieu \cite{FeVi06}, for several case studies).
Similarly to traditional regression problems, model (\ref
{funclinearregression}) may arise under different experimental designs.
We assume a random design of the explanatory curves, where
$X_{1},\ldots, X_{n}$ is a sequence of identically distributed random
functions with the same
distribution as a generic $X$. The main assumption on $X$ is that it
is a second-order variable, that is, $\mathbb{E} (\int_I
X^2(t)\,dt ) < + \infty,$ and it is assumed moreover that $\mathbb
{E}(X_i(t)\varepsilon_i)=0$
for almost every $t\in I$.
This situation has been considered, for instance, in Cardot, Ferraty and Sarda
\cite{CaFeSa03} and M\"uller and Stadtm\"uller \cite{MuSt05} for
independent variables, while correlated functional variables are
studied in
Bosq \cite{Bo00}. Our analysis is based on a
general framework without any assumption of independence of the
$X_i$'s. We will, however, assume independence between the $X_i$'s and
the $\varepsilon_i$'s in our theoretical results
in Sections \ref{sec3} and \ref{sec4}.

The main problem in functional linear regression is
to derive an estimator $\widehat{\alpha}$ of the unknown slope function
$\alpha$. However, estimation of $\alpha$ in
(\ref{funclinearregression}) belongs to the class of ill-posed inverse
problems. Writing (\ref{funclinearregression}) for generic
variables $X$,
$Y$ and $\varepsilon$, multiplying both sides by $X-\mathbb{E}(X)$
and then taking expectations leads to
%
\begin{eqnarray}\label{normalequation}
&&\mathbb{E}\bigl(\bigl(Y-\mathbb{E}(Y)\bigr)\bigl(X-\mathbb{E}(X)\bigr)\bigr)\nonumber\\[-8pt]\\[-8pt]
&&\qquad = \mathbb{E}
\biggl(\int_{I} \alpha(t) \bigl(X(t)-\mathbb{E}(X)(t)\bigr) \,dt
\bigl(X-\mathbb{E}(X)\bigr)\biggr)=: \Gamma(\alpha).\nonumber
\end{eqnarray}
The normal equation (\ref{normalequation}) is the continuous equivalent
of normal equations in the multivariate linear model.
Estimation of $\alpha$ is thus linked with the inversion of the
covariance operator $\Gamma$ of $X$ defined in (\ref{normalequation}).
But, unlike the finite dimensional case, a bounded inverse for $\Gamma
$ does not exist since it is a compact linear operator
defined
on the infinite dimensional space $L^2(I)$. This corresponds to the
setup of ill-posed inverse problems (with the additional difficulty that
$\Gamma$ is unknown). As a consequence, the parameter $\alpha$ in
(\ref{funclinearregression}) is not identifiable without additional
constraint. Actually, a necessary and sufficient condition
under which a unique solution for
(\ref{funclinearregression})--(\ref{normalequation}) exists in the
orthogonal space of $\operatorname{ker}(\Gamma)$
and is given by
$\sum_{r}(\frac{\mathbb{E}((Y-\mathbb{E}(Y))\int_I(X(t)-\mathbb
{E}(X)(t))\zeta_r(t)\,dt)}{\lambda_r})^2 < +\infty,$ %
where $(\lambda_r,\zeta_r)_r$ are the eigenelements of $\Gamma$
(see Cardot, Ferraty and Sarda \cite{CaFeSa03} or He, M\"uller and
Wang \cite{HeMuWa00} for a functional response).
The set of solutions is the set of functions $\alpha$ which can be
decomposed as a sum of the unique element of the orthogonal space of
$\operatorname{ker}(\Gamma)$ satisfying (\ref{normalequation}) and any
element of $\operatorname{ker}(\Gamma)$.

It follows from these arguments that any sensible procedure for
estimating $\alpha$ (or, more precisely, of its identifiable part)
has to involve regularization procedures.
Several authors have proposed estimation procedures where
regularization is obtained in two main ways. The first one is based
on the Karhunen--Lo\`eve expansion of $X$ and leads to regression on
functional principal components:
see Bosq~\cite{Bo00}, Cardot, Mas and Sarda \cite{CaMaSa07} or M\"
uller and Stadtm\"uller \cite{MuSt05}. It consists in projecting the
observations on a finite dimensional space spanned by eigenfunctions
of the (empirical) covariance operator $\Gamma_n$.
For the second method, regularization is obtained through a penalized
least squares approach after expanding $\alpha$
in some basis (such as splines): see Ramsay and Dalzell \cite{RaDa91},
Eilers and Marx \cite{EiMa96}, Cardot, Ferraty and
Sarda \cite{CaFeSa03} or Li and Hsing \cite{Li06}.
We propose here to use a smoothing splines approach prolonging a
previous work from Cardot et al. \cite{CaCrKnSa07}.

Our estimator is described in Section \ref{sec2}. Note that (\ref
{funclinearregression}) implies
that $Y_i-\overline{Y}=\int_I \alpha(t)[X_i(t)-\overline
{X}(t)]\,dt+\varepsilon_i-\bar\varepsilon$.
Based on the observation times $t_{1} < \cdots < t_{p}$, we rely on minimizing
the residual sum of squares $\sum_i ( Y_i -\overline{Y}
- \frac{1}{p}\sum_{j=1}^p a(t_j)(X_i(t_j)-\overline{X}(t_j)
))^2$ subject to a roughness penalty. A slight modification
of the usual penalty term is applied in order to guarantee the
existence of the estimator under general conditions.
The proposed estimator $\widehat{\alpha}$ is then a natural spline with
knots at the observation points $t_j$. An\vspace*{2pt} estimator of the intercept
$\alpha_0=\mathbb{E}(Y)-\int_I \alpha(t)\mathbb{E}( X)(t)]\,dt$ is
given by $\widehat{\alpha_0}=\overline{Y}-\int_I \widehat{\alpha
}(t)\overline{X}(t)\,dt$.
For simplicity, we will
assume that $t_{1} < \cdots < t_{p}$ are equispaced, but the
methodology can easily be generalized to other situations.
It must be emphasized, however, that our study does not cover the case
of sparse points for which other techniques have to be
envisaged; for this specific problem, see the work from
Yao, M\"uller and Wang~\cite{YaMuWa05}.

In Section \ref{sec3} we present a detailed asymptotic theory of the behavior of
our estimator for large values of $n$ and $p$.
The distance between $\widehat{\alpha}$ and $\alpha$ is evaluated
with respect to $L^2$ semi-norms induced by the operator $\Gamma$,
$\| u\|^2_\Gamma=\langle\Gamma u,u\rangle$ with
$\langle u,v\rangle=\int_I u(t)v(t)\,dt$, or its discretized or
empirical versions (see, e.g.,
Cardot, Ferraty and Sarda \cite{CaFeSa03} or M\"uller and Stadtm\"uller \cite
{MuSt05} for similar setups).
By using these semi-norms we explicitly concentrate on analyzing the
estimation error only for the identifiable part of the structure of
$\alpha$ which is relevant for prediction.
Indeed, it will be shown in Section \ref{sec3} that $\|\widehat{\alpha}-\alpha\|
_\Gamma^2$ determines the rate of convergence of
the error in predicting the conditional mean $\alpha_0+\int_I \alpha
(t)X_{n+1}(t)\,dt$ of $Y_{n+1}$
for any new random function $X_{n+1}$ possessing the same distribution
as $X$ and independent of
$X_1,\ldots,X_n$:
%
\begin{eqnarray}\label{prederror}
&&\mathbb{E} \biggl(\biggl(\widehat{\alpha_0}+\int_I \widehat{\alpha
}(t)X_{n+1}(t)\,dt-
\alpha_0-\int_I \alpha(t)X_{n+1}(t)\,dt\biggr)^2\Big| \widehat{\alpha_0},\widehat{\alpha}
\biggr)\nonumber\\[-8pt]\\[-8pt]
&&\qquad= \|\widehat{\alpha}-\alpha\|_\Gamma^2+O_P(n^{-1}).\nonumber
\end{eqnarray}

We first derived optimal rates of convergence with respect to the $L^2$
semi-norms induced by $\Gamma$ in a quite general setting which substantially
improved existing results in the literature as well as bounds
obtained for this estimator in a previous paper (see Cardot
et al. \cite{CaCrKnSa07}). If $\alpha$ is $m$-times
continuously differentiable, then it is shown that rates of
convergence for our estimator are of
order $n^{-(2m+2q+1)/(2m+2q+2)}$, where
the value of $q>0$ depends on the structure of the distribution of
$X$. More precisely,
$q$ quantifies the rate of decrease $\sum_{r=k+1}^\infty\lambda
_r=O(k^{-2q})$ as
$k\rightarrow\infty$, where $\lambda_1\ge\lambda_2\ge\cdots$
are the eigenvalues of the covariance operator $\Gamma$. If, for
example, $X$ is a.s. twice continuously differentiable,
then $q\ge2$.
As a second step, we show that these rates of convergence are
optimal in the sense that they are minimax over large classes of
distributions of~$X$ and of functions
$\alpha$. No alternative estimator can globally achieve faster rates
of convergence in these classes.

In an interesting paper Cai and Hall \cite{CaHa06} derive rates of
convergence on the error $\alpha_0+\langle\alpha, x\rangle- \widehat
{\alpha}_0-\langle\widehat\alpha, x\rangle$ for a
pre-specified, fixed function $x$. Their approach is based on
regression with respect to functional
principal components and the derived rates are shown to be optimal
with respect to this methodology.
At first glance this setup seem to be close, but
due to the fact that explanatory variables are of infinite dimension, inference
on \textit{fixed} functions $x$ cannot generally be used to derive
optimal rates of convergence of the
prediction error (\ref{prederror}) for \textit{random} functions
$X_{n+1}$. We also want to emphasize that
in the present paper we do not
consider the convergence of $\hat\alpha$ with respect to the usual
$L^2$ norm. Analyzing $\|\widehat{\alpha}-
\alpha\|^2=\int_I(\widehat{\alpha}(t)-\alpha(t))^2\,dt$
instead of
$\|\widehat{\alpha}-
\alpha\|_\Gamma^2$ must be seen statistically
as a very different problem, and under our general assumptions it only
follows that
$\|\widehat{\alpha}-
\alpha\|^2$ is
bounded in probability
(see the proof of Theorem \ref{thm2}). It appears that to get stronger results
one needs additional
conditions linking the ``smoothness'' of $\alpha$ and of the curves
$X_i$ as derived in a
recent work by Hall and Horowitz \cite{HaHo07}. A detailed discussion
of these issues is given in Section \ref{sec32}.

In practice the functional values $X_i(t_j)$ are often not directly
observed; there exist only noisy observations $W_{ij}=X_i(t_j)+\delta_{ij}$
contaminated with random errors~$\delta_{ij}$. In Section \ref{sec4}, we
consider a modified functional linear model adapting to such
situations. In this
errors-in-variable context, we use a corrected estimator as introduced
in Cardot et al. \cite{CaCrKnSa07} which can be seen as a modified
version of the so-called \textit{total least squares} method for
functional data. We show again the good asymptotic performance of the method
for a sufficiently dense grid of discretization points.

We devote Section \ref{sec5} to the application of the proposed
estimation procedure to the prediction of
the peak of pollution from the curve of pollutant indicators collected
the preceding day. Finally, the
proofs of our results can be found in Section~\ref{sec6}.

\section{Smoothing splines estimation of the functional coefficient}
\label{sec2}

As explained in the \hyperref[sec1]{Introduction}, we will assume that the functions
$X_i$ are observed at $p$
equidistant points $t_1,\ldots,t_p\in I$. In order to simplify further
developments, we will take $I=[0,1]$
so that $t_1=\frac{1}{2p}$ and $t_{j} - t_{j-1} = \frac{1}{p}$
for all $j = 2, \ldots, p$.

Our estimator of $\alpha$ in (\ref{funclinearregression}) is a
generalization
of the well-known smoothing splines estimator in univariate
nonparametric regression.
It relies on the implicit assumption that the underlying function
$\alpha$ is sufficiently smooth
as, for example, $m$-times continuously differentiable ($m=1,2,3,\ldots$).

For any smooth function $a$ the discrete sum $\frac{1}{p}\sum_{j=1}^p
a(t_j)X_i(t_j)$
is used to approximate the integral $\int_0^1 a(t)X_i(t)\,dt$ in (\ref
{funclinearregression}), whereas expectations are estimated by the
sample means $\overline{Y}$ and $\overline{X}$,
and an estimate is obtained by minimizing the sum of squared residuals
$( Y_i-\overline{Y}- \frac{1}{p}\sum_{j=1}^p a(t_j)(X_i(t_j)-\overline
{X}(t_j)))^2$ subject to a roughness penalty.
More precisely,
for some $m=1,2,\ldots$ and a smoothing parameter $\rho>0$,
an estimate $\widehat{\alpha}$ is determined by minimizing
%
\begin{eqnarray}\label{splinedef}
&&
\frac{1}{n}\sum_{i=1}^n \Biggl( Y_i-\overline{Y}- \frac{1}{p}\sum
_{j=1}^p
a(t_j)\bigl(X_i(t_j)-\overline{X}(t_j)\bigr)\Biggr)^2\nonumber\\[-8pt]\\[-8pt]
&&\qquad{}
+\rho\Biggl(\frac{1}{p}\sum_{j=1}^p \pi_a^2(t_j) +\int_{0}^{1}
\bigl(a^{(m)}(t)\bigr)^2\,dt\Biggr)\nonumber
\end{eqnarray}
over all functions $a$ in the Sobolev space $W^{m,2}([0,1])\subset
L^2([0,1])$, where $\pi_a(t)=\sum_{l=1}^m
\beta_{a,l}t^{l-1}$ with $\sum_{j=1}^p (a(t_j)-\pi
_a(t_j))^2=\min_{\beta_1,\ldots,\beta_m}
\sum_{j=1}^p (a(t_j)-\sum_{l=1}^m
\beta_{l}t^{l-1})^2$.

Obviously, $\pi_a$ denotes the best possible approximation of
$(a(t_1),\ldots,a(t_p))$ by a polynomial of
degree $m-1$. The extra term $\frac{1}{p}\sum_{j=1}^p \pi_a(t_j)^2$ in
the roughness penalty is unusual
and does not appear in traditional smoothing splines approaches. It
will, however, be shown below that
this term is necessary to guarantee existence of a unique solution in a
general context
without any additional
assumptions on the curves $X_i$.

It is quite easily seen that any solution $\widehat{\alpha}$ of (\ref
{splinedef})
has to be an element of the space $NS^{m} (t_{1},\ldots, t_{p})$ of
\textit{natural splines} of order
$2m$ with knots at $t_{1},\ldots, t_{p}$. Recall that $NS^{m} (t_{1},\ldots, t_{p})$
is a $p$-dimensional linear space of
functions with $v^{(m)} \in L^{2}([0,1])$ for any $v \in NS^{m}
(t_{1},\ldots, t_{p})$.
Let $\mathbf{b}(t) = ( b_{1}(t),\ldots, b_{p} (t) )^{\tau}$
be a functional basis
of $NS^{m} (t_{1},\ldots, t_{p})$. A discussion of
several possible basis function expansions
can be found in Eubank \cite{Eu88}. An important property of natural
splines is that
there exists a canonical one-to-one mapping between $\mathbb{R}^p$ and
the space $NS^{m}
(t_{1},\ldots, t_{p})$ in the following way: for any vector
$\mathbf{w} = ( w_{1}, \ldots, w_{p} )^{\tau} \in
\mathbb{R}^{p}$,
there exists a unique natural spline interpolant $s_{\mathbf{w}}$ with
$s_{\mathbf{w}} (t_{j})= w_{j}$, $j=1,\ldots, p$. With $\mathbf{B}$ denoting
the $p \times p$ matrix with elements $b_{i}(t_{j})$, $s_{\mathbf
{w}}$ is
given by
%
\begin{equation}\label{natspline}
s_{\mathbf{w}} (t) = \mathbf{b}(t)^{\tau} ( \mathbf{B}^{\tau}
\mathbf{B}
)^{-1} \mathbf{B}^{\tau} \mathbf{w}.
\end{equation}
The important property of such a spline interpolant is the
fact that
\begin{eqnarray}\label{natsmooth}
&&\int_{0}^1 s_{\mathbf{w}}^{(m)}(t)^{2} \,dt \leq\int_{0}^1
f^{(m)}(t)^{2} \,dt\\
&&\eqntext{\mbox{for any other function } f\in W^{m,2}([0,1])}
\\[-4pt]
&&\eqntext{\mbox{ with } f(t_{j}) = w_{j}, j=1,\ldots, p.}
\end{eqnarray}

Note that in (\ref{splinedef}) only the integral $\int_{0}^{1}
a^{(m)}(t)^2\,dt$ depends
on the values of $a$ in the open intervals $(t_{j-1},t_j)$ between grid
points. It therefore
follows from (\ref{natsmooth}) that $\widehat{\alpha}= s_{\widehat
{\bolds{\alpha}}}$, where
$\widehat{\bolds{\alpha}}=(\widehat{\alpha}(t_1),\ldots,\widehat
{\alpha}(t_p))^\tau\in\mathbb{R}^{p}$
minimizes
%
\begin{eqnarray}\label{splinedef1}
&&\frac{1}{n}\sum_{i=1}^n \Biggl( Y_i-\overline{Y}- \frac{1}{p}\sum
_{j=1}^p
a(t_j)\bigl(X_i(t_j)-\overline{X}(t_j)\bigr)\Biggr)^2\nonumber\\[-8pt]\\[-8pt]
&&\qquad{} +\rho\Biggl(\frac{1}{p}\sum_{j=1}^p \pi_a^2(t_j) +\int_{0}^{1}
s_{\mathbf{a}}^{(m)}(t)^2\,dt\Biggr);\nonumber
\end{eqnarray}
with respect to all vectors $\mathbf{a}=(a(t_1),\ldots,a(t_p))^\tau
\in\mathbb{R}^{p}$.

A closer study of $\widehat{\bolds{\alpha}}$ requires the use of
matrix notation:
$\mathbf{Y} = ( Y_{1}-\overline{Y}, \ldots,
Y_{n}-\overline{Y} )^{\tau}$,
$\mathbf{X}_{i} = ( X_{i} (t_{1})-\overline{X}(t_1), \ldots,
X_{i} (t_{p})-\overline{X}(t_p) )^{\tau}$ for
all $i=1, \ldots, n$, $\bolds{\alpha} = ( \alpha(t_{1}),
\ldots,
\alpha(t_{p}) )^{\tau}$, $\bolds{\varepsilon} =
( \varepsilon_{1}-\overline{\varepsilon}, \ldots, \varepsilon_{n}-\overline
{\varepsilon} )^{\tau}$ and let
$\mathbf{X}$ be the $n \times p$ matrix with a general term
$X_{i} (t_{j})-\overline{X}(t_j)$ for all $i=1, \ldots, n$,\break $j=1,
\ldots, p$.
Moreover, $\mathbf{P}_m$ will denote the $p \times p$ projection\break
matrix projecting into the
$m$-dimensional linear space $E_{m}:=\{\mathbf{w}=(w_1,\ldots
,\break w_p)^\tau\in\mathbb{R}^{p}|
 w_{j} = \sum_{l=1}^{m}\theta_{l} t_{j}^{l-1}$, $j=1,\ldots, p\}$ of
all (discretized) polynomials of
degree $m-1$. By (\ref{natspline}), we have $\int_{0}^1 s_{\mathbf
{a}}^{(m)}(t)^{2} \,dt =
\mathbf{a}^\tau\mathbf{A}_{m}^{*} \mathbf{a}$,
where
$\mathbf{A}_{m}^{*} = \mathbf{B} ( \mathbf{B}^{\tau} \mathbf{B}
)^{-1}
[ \int_{0}^1 \mathbf{b}^{(m)}(t) \mathbf{b}^{(m)}(t)^{\tau} \,dt ]
( \mathbf{B}^{\tau} \mathbf{B} )^{-1} \mathbf{B}^{\tau}$ is a
$p\times p$ matrix.

When defining $\mathbf{A}_{m}:=\mathbf{P}_m+p\mathbf{A}_{m}^{*}$,
minimizing (\ref{splinedef1}) is equivalent to solving
%
\begin{equation}\label{nonnoisypbmino}
\min_{\mathbf{a}\in\mathbb{R}^{p}}\biggl\{ \frac{1}{n} \biggl\|
\mathbf{Y}
- \frac{1}{p} \mathbf{X} \mathbf{a} \biggr\|^{2} +
\frac{\rho}{p} \mathbf{a}^\tau\mathbf{A}_{m}\mathbf{a}
\biggr\} ,
\end{equation}
where $\|\cdot\|$ stands for the usual Euclidean norm. The
solution is given by
%
\begin{equation}\label{solX}\qquad
\widehat{\bolds{\alpha}} = \frac{1}{n p} \biggl( \frac{1}{n p^{2}}
\mathbf{X}^{\tau} \mathbf{X} + \frac{\rho}{p} \mathbf{A}_{m} \biggr)^{-1}
\mathbf{X}^{\tau} \mathbf{Y} = \frac{1}{n} \biggl( \frac{1}{n p} \mathbf
{X}^{\tau}
\mathbf{X} + \rho\mathbf{A}_{m} \biggr)^{-1} \mathbf{X}^{\tau} \mathbf
{Y} .
\end{equation}
Then $\widehat{\alpha}= s_{\widehat{\bolds{\alpha}}}$ constitutes
our final estimator of $\alpha$ while
$\widehat{\alpha_0}=\overline{Y}-\langle\widehat{\alpha},\overline
{X}\rangle$ is used to estimate the intercept $\alpha_0$.
Based on a somewhat different development, this estimator of $\alpha$
has already been proposed by
Cardot et al. \cite{CaCrKnSa07}.

In order to verify existence of $\widehat{\bolds{\alpha}}$,
let us first cite some properties of the eigenvalues of $p \mathbf
{A}_{m}^{*}$ which have been studied by many authors
(see Eubank \cite{Eu88}). For instance, in Utreras \cite{Ut83}, it is
shown that this matrix has exactly $m$ zero
eigenvalues $\mu_{1,p} = \cdots = \mu_{m,p} = 0$. The corresponding
$m$-dimensional eigenspace is the space
$E_m$ of discretized polynomials as defined above.
The $p-m$ nonzero eigenvalues
$0 < \mu_{m+1,p} < \cdots < \mu_{p,p}$ are such that there exist
constants $0<D_0<D_1<\infty$ such that
$D_0 \leq\mu_{j+m,p} (\pi j)^{-2m} \leq D_1$ for $j = 1, \ldots,
p-m$ and all sufficiently large $p$. Therefore, there exist some
constant $0<C_0<+\infty$ and some $p_0\in\{0,1,2,\ldots\}$ such that
for all $p\ge p_0$ and $k=0,\ldots,p-m-1$
\begin{eqnarray}\label{egams}
k^{2m} \frac{1}{\mu_{k+m+1,p}} \leq C_0.
\end{eqnarray}
We can conclude that all eigenvalues of the matrix $\mathbf
{A}_{m}$ are
strictly positive, and existence as well as uniqueness of the solution
(\ref{solX}) of the minimization problem~(\ref{nonnoisypbmino}) are straightforward consequences. Note that
\hyperref[sec1]{Introduction} of the additional term
$\frac{1}{p}\sum_{j=1}^p \pi_a(t_j)^2$ in (\ref{splinedef}) is crucial.
Dropping this term in
(\ref{splinedef}) as well as (\ref{splinedef1}) results in replacing
$\mathbf{A}_{m}$ by $p\mathbf{A}_{m}^{*}$
in (\ref{nonnoisypbmino}). Existence of a solution then cannot be
guaranteed in a general context since, due
to the $m$ zero eigenvalues of $p\mathbf{A}_{m}^{*}$, the matrix
$( \frac{1}{n p^{2}}
\mathbf{X}^{\tau} \mathbf{X} + \rho\mathbf{A}_{m}^{*} )$ may not be invertible.

\begin{remark*}
Our requirement of equidistant grid points $t_j$ has to be seen as a
restrictive condition. There are many applications where
the functions $X_i$ are only observed at varying numbers $p_i$ of
irregularly spaced points $t_{i1}\le\cdots\le t_{ip_i}$. Then our
estimation procedure is not directly applicable. Fortunately there
exists a fairly simple modification.
Define a smooth function $\widetilde X_i\in L^2([0,1])$ by smoothly
interpolating the observations (e.g., using natural splines) such that
$\widetilde X_i(t_{ij})=X_i(t_{ij})$, $j=1,\ldots,p_i$. Then define
$p>\max\{p_1,\ldots,p_n\}$ equidistant grid points $t_1,\ldots,t_p$, and
determine an estimator $\widehat{\alpha}$ by applying the smoothing
spline procedure (\ref{splinedef}) with
$\frac{1}{p}\sum_{j=1}^p a(t_j)(X_i(t_j)-\overline{X}(t_j))$ being
replaced by $\frac{1}{p}\sum_{j=1}^p a(t_j)(\widetilde
X_i(t_j)-\overline{\widetilde{X}}(t_j))$.
For example, in the case of a random design with i.i.d. observations
$t_{ij}$ from a strictly positive design density on $I$,
it may be shown that the asymptotic results of Section \ref{sec3} generalize to
this situation if $\min\{p_1,\ldots,$ $p_n\}$ is sufficiently large
compared to $n$. A detailed analysis is not in the scope
of the present paper.
\end{remark*}

\section{Theoretical results} \label{sec3}

\subsection{Rates of convergence for smoothing splines estimators}
\label{sec31}

We will denote the standard inner product of the Hilbert space
$L^2([0,1])$ by
$\langle f,g\rangle=\int_0^1f(t)g(t)\,dt$ and $\|\cdot\|$ by its associated
norm. As outlined in the \hyperref[sec1]{Introduction}, our analysis is based on evaluating
the error between $\widehat{\alpha}$ and $\alpha$ with respect to the semi-norm $\|\cdot\|_\Gamma$
defined in Section \ref{sec1},
\[
\| u \|_{\Gamma}^{2} := \langle\Gamma u,u\rangle,\qquad u\in
L^2([0,1]),
\]
where $\Gamma$ is the \textit{covariance operator} of $X$
given by
\[
\Gamma u := \mathbb{E} \bigl( \bigl\langle\bigl(X -\mathbb{E}(X)\bigr), u \bigr\rangle
\bigl(X-\mathbb{E}(X)\bigr) \bigr), \qquad u\in L^{2}([0,1]).
\]

The above $L^2$ semi-norm has already been
used in similar contexts as the one studied in the present paper;
see, for example, Wahba \cite{Wa77}, Cardot, Ferraty and Sarda \cite
{CaFeSa03} or M\"uller and Stadtm\"uller \cite{MuSt05}.
By (\ref{prederror}) the asymptotic behavior of $\| \widehat
{\alpha}-\alpha\|_{\Gamma}^{2}$
constitutes a major object of interest, since it quantifies the
leading term in the expected squared prediction error
for a new random function $X_{n+1}$.

As first steps, we will consider in Theorems \ref{thm1} and \ref{thm2}
the error between $\widehat{\alpha}$ and $\alpha$ with respect to
simplified versions of the above semi-norm: the
discretized empirical semi-norm defined for any $\mathbf{u}\in
\mathbb{R}^p$ as
\[
\| \mathbf{u} \|_{\Gamma_{n,p}}^{2} := \frac{1}{p}
\mathbf{u}^{\tau}
\biggl( \frac{1}{n p} \mathbf{X}^{\tau} \mathbf{X} \biggr) \mathbf
{u},
\]
and the empirical semi-norm defined for any $u \in
L^{2}([0,1])$ as
\[
\| u \|_{\Gamma_{n}}^{2} := \frac{1}{n} \sum_{i=1}^{n}
\langle(X_{i}-\overline{X}) , u \rangle^{2} = \langle\Gamma_{n} u , u
\rangle,
\]
where $\Gamma_n$ is the empirical covariance operator from $X_1,\ldots
,X_n$ given by
\[
\Gamma_{n} u := \frac{1}{n} \sum_{i=1}^{n} \langle(X_{i}-\overline{X})
, u \rangle(X_{i}-\overline{X}).
\]

Obviously, $\|\widehat{\bolds{\alpha}}-\bolds{\alpha}
\|_{\Gamma_{n,p}}^{2}
=\frac{1}{n}\sum_i [\frac{1}{p} \sum_{j=1}^p(\widehat{\alpha
}(t_j)-\alpha(t_j))(X_i(t_j)-\overline{X}(t_j))]^2$
and
$\| \widehat{\alpha}-\alpha\|_{\Gamma_{n}}^{2} =\frac
{1}{n}\sum_i [\int_I(\widehat{\alpha}(t)-\alpha(t))(X_i(t)-\overline
{X}(t))\,dt]^2$
quantify different modes of convergence of $\langle\widehat{\alpha},
X-\overline{X}\rangle$ to $\langle\alpha,(X-\overline{X})\rangle$.

As mentioned in Section \ref{sec2}, the function $\alpha$ is required to have a
certain degree of regularity.
Namely, it satisfies the following assumption for some $m\in\{1,2,\ldots
\}$:
\renewcommand{\theequation}{A.1}
\begin{equation}\label{assa1}
\alpha\mbox{ is $m$-times differentiable and $\alpha^{(m)}$ belongs to
$L^{2}([0,1])$.}
\end{equation}
Let $C_{1} = \int_0^1 \alpha^{(m)}(t)^{2} \,dt$ and $C_{2}^{*}
= \int_0^1 \alpha(t)^{2} \,dt$. By
construction of $\mathbf{P}_{m}$, $\mathbf{P}_{m} \bolds{\alpha}$
provides the best approximation (in a least squares sense)
of $\bolds{\alpha}$ by (discretized) polynomials of degree $m-1$,
and $\frac{1}{p} \bolds{\alpha}^{\tau} \mathbf{P}_{m}
\bolds{\alpha} \leq\frac{1}{p}\bolds{\alpha}^{\tau}\mathbf{A}_{m}
\bolds{\alpha} \longrightarrow C_{2}^{*}$ as $p \rightarrow\infty
$. Let $C_{2}$ denote
an arbitrary constant with $C_{2}^{*} < C_{2} < \infty$. There then
exists a $p_{1} \in\{0,1,\ldots\}$ with $p_{1}
\geq p_{0}$ such
that $\frac{1}{p}\bolds{\alpha}^{\tau} \mathbf{P}_{m} \bolds
{\alpha} \leq C_{2}$ for all $p \geq p_{1}$.

Recall that our basic setup implies that $X_{1},\ldots, X_{n}$ are
identically distributed random functions with the same
distribution as a generic variable $X$.
Expected values $\mathbb{E}_\varepsilon(\cdot)$ as stated in the theorems
below will refer to the probability
distribution induced by the random variable $\varepsilon$, that is,
they stand for conditional expectation
given $X_{1},\ldots, X_{n}$. We assume moreover that $\varepsilon_i$ is
independent of the~$X_i$'s. In the following, for any real positive number
$x$, $[x]$ will denote the smallest integer which is larger than $x$.
In addition,
let $\lambda_{x,1}\ge\lambda_{x,2}\ge\cdots\ge\lambda_{x,p}\ge0$
denote the eigenvalues of the matrix $\frac{1}{n p} \mathbf{X}^{\tau} \mathbf{X}$.
We start with a theorem giving finite sample bounds for bias and
variance of the estimator $\widehat{\bolds{\alpha}}$ with respect
to the semi-norm~$\| \cdot \|_{\Gamma_{n,p}}$.

\begin{theorem} \label{thm1}
Under assumption \textup{(A.1)} and the above definitions of $C_{0}$, $C_{1}$,
$C_{2}$, $p_{1}$, the following bounds hold
for all $n=0,1,\ldots,$ all $p \geq p_{1}$, all $\rho>n^{-2m}$ and
every $n\times p$ matrix $\mathbf{X}=(X_i(t_j))_{i,j}$:
%
\setcounter{equation}{0}
\renewcommand{\theequation}{\arabic{section}.\arabic{equation}}
\begin{eqnarray}\label{resthm1}
\| \mathbb{E}_{\varepsilon} ( \widehat{\bolds{\alpha}}) -
\bolds{\alpha} \|_{\Gamma_{n,p}}^{2} & \leq
& 2\rho\biggl( \frac{1}{p} \bolds{\alpha}^{\tau} \mathbf{P}_{m}
\bolds{\alpha} + C_{1} \biggr)
+\frac{4}{n}\sum_{i=1}^n(d_i-\overline{d})^2 \nonumber\\[-8pt]\\[-8pt]
& \leq& \rho( C_{2} + C_{1} )+\frac{4}{n}\sum
_{i=1}^n(d_i-\overline{d})^2 ,\nonumber
\end{eqnarray}
as well as
%
\begin{equation}\label{resthm12}
\mathbb{E}_{\varepsilon} \bigl(\| \widehat{\bolds{\alpha}} -
\mathbb{E}_{\varepsilon} ( \widehat{\bolds{\alpha}}
)\|_{\Gamma_{n,p}}^{2} \bigr) \leq\frac{\sigma^{2}_{\varepsilon
}}{n} \Bigl( m+\bigl[\rho^{-{1}/({2m+2q+1})}\bigr]
(2+C\cdot C_0)\Bigr),
\end{equation}
for any $C>0$ and $q\ge0$ with the property that $ \sum
_{j=k+1}^p\lambda_{x,j}\le C\cdot k^{-2q}$ holds for $k:=[\rho^{-
{1}/({2m+2q+1})}]$.
\end{theorem}

The rate of convergence
of $\| \widehat{\bolds{\alpha}} - \bolds{\alpha} \|_{\Gamma
_{n,p}}^{2}$ thus
depends on assumptions on the distribution of $X$ and on the size of
the discretization error.
In order to complement our basic setup,
we will rely on the following conditions:

\begin{enumerate}[(A.2)]
\item[(A.2)] There exists some constant $\kappa$, $0<\kappa<1,$ such
that for every $\delta>0$, there exists a constant $C_{3} < + \infty$
such that
\[
\mathbb{P}\bigl(|X(t) - X(s)| \leq C_{3} |t-s|^{ \kappa}, t,s\in
I\bigr)\ge1-\delta.
\]
\item[(A.3)] For some constant $C_{4}<\infty$ and
all $k=1,2,\ldots$ there is a $k$-dimensional linear subspace $\mathcal{L}_k$ of $L^2([0,1])$ with
\[
\mathbb{E}\biggl( \inf_{f\in\mathcal{L}_k} \sup_t |X(t) -f(t)|^2
\biggr)\leq C_{4} k^{-2q}.
\]
\end{enumerate}

Before proceeding any further, let us consider assumption (A.3) more
closely. The following lemma provides a link
between assumption (A.3) and the degree of smoothness of the random
functions $X_i$.

\begin{lemma} \label{lem1}
For some $q_1=0,1,2,\ldots$ and $0\le r_2\le1$ assume that $X$
is almost surely $q_1$-times continuously differentiable and that
there exists some $C_5<\infty$ such that
\[
\mathbb{E}\biggl( \sup_{|t-s|\le d} \bigl|X^{(q_1)}(t) -
X^{(q_1)}(s)\bigr|^2\biggr) \leq C_{5} d^{ 2r_2}
\]
holds for all $d>0$. There then exists a constant $C_6<\infty$,
depending only on $q_1$, such that
for all $k=1,2,\ldots$
\[
\mathbb{E}\biggl( \inf_{f\in\mathcal{E}_k} \sup_t |X(t) -f(t)|^2\biggr)
\leq C_{6}C_{5} k^{-2(q_1+r_2)},
\]
where $\mathcal{E}_k$ denotes the space of all polynomials of order $k$ on $[0,1]$.
\end{lemma}

\begin{pf}
The well-known Jackson's inequality in approximation
theory implies the existence of some $C_6<\infty$,
only depending on $q_1$, such that for all $k=1,2,\ldots$
\begin{eqnarray*}
\inf_{f\in\mathcal{E}_k}\sum_{j=1}^p \bigl(X(t_j)-f(t_j)\bigr)^2
\le C_6k^{-2q_1} \sup_{|t-s|\le1/k} \bigl|X^{(q_1)}(t) - X^{(q_1)}(s)\bigr|^2
\end{eqnarray*}
holds with probability 1. The lemma is an immediate consequence.
\end{pf}

The lemma implies that if assumption (A.2) can be replaced by the
stronger requirement
$\mathbb{E}( \sup_{|t-s|\le d} |X(t) - X(s)|^2) \leq C_{5} d^{ -2r_2}$, $d>0$,
then assumption (A.3) necessarily holds for some $q\ge\kappa$.
Indeed, $q\gg\kappa$ will result
from a very high degree of smoothness
of $X_i$.

On the other hand, assumption (A.3) only requires that the functions
$X_i$ be well approximated by
\textit{some} arbitrary low dimensional linear function spaces (not
necessarily polynomials). Even if $X_i$ are \textit{not} smooth,
assumption (A.3) may be satisfied for a large value of $q$ (the
Brownian motion provides an example).

Theorem \ref{thm1} together with assumptions (A.2) and (A.3) now allows us to
derive rates of convergence of our estimator $\widehat{\alpha}$.
First note that
assumption (A.3) determines the rate of decrease of the eigenvalues
$\lambda_{x,j}$ of $\frac{1}{n p} \mathbf{X}^{\tau} \mathbf{X}$.
For\vspace*{-1pt} any \mbox{$k$-dimensional} linear space $\mathcal{L}_k\subset L^2([0,1])$,
let $\mathcal{P}_k$ denote the corresponding $p\times p$ projection
matrix projecting
into the $k$-dimensional subspace $\mathcal{L}_{k,p}=\{ v\in\mathbb{R}^p|
v=(f(t_1),\ldots,f(t_p))^\tau,  f\in\mathcal{L}_k\}$.
Basic properties of eigenvalues and eigenvectors then imply that
\setcounter{equation}{2}
\renewcommand{\theequation}{\arabic{section}.\arabic{equation}}
\begin{eqnarray}\label{H31}
&&\sum_{j=k+1}^p \lambda_{x,j}\le \inf_{ \mathcal{P}_k}\operatorname{Tr}
\biggl(
(\mathbf{I}_p-\mathcal{P}_k)\frac{1}{n p} \mathbf{X}^{\tau} \mathbf{X}
\biggr)\nonumber\\[-8pt]\\[-8pt]
&&\qquad=\frac{1}{n p}\sum_{i=1}^n \inf_{f\in\mathcal{L}_k}\sum_{j=1}^p
\bigl(X_i(t_j)-\overline{X}-f(t_j)\bigr)^2,\nonumber
\end{eqnarray}
and assumption (A.3) implies that for any
$\delta>0$ there exists a $C_\delta<\infty$ such that $P(\sum
_{j=k+1}^p\lambda_{x,j}\le C_{\delta} k^{-2q})\ge1-\delta$.

Assumptions (A.1) and (A.2) obviously lead to
%
\begin{equation} \label{H3a}
\frac{1}{n} \sum_{i=1}^{n} (d_{i}-\overline{d})^{2} =
O_P ( p^{-2 \kappa} ).
\end{equation}
If $n,p\rightarrow\infty$, $\rho\rightarrow0$, $1/(n\rho
)\rightarrow0$, then
relations (\ref{resthm1}), (\ref{resthm12}) and (\ref{H31}) imply that
\[
\| \widehat{\bolds{\alpha}} - \bolds{\alpha} \|
_{\Gamma_{n,p}}^{2} =
O_{P} \bigl( \rho+ \bigl(n\rho^{{1}/({2m+2q+1})}\bigr)^{-1}+p^{-2 \kappa}
\bigr).
\]
In the following we will require that $p$ is sufficiently large
compared to $n$ so that the discretization error is negligible.
It therefore suffices that $n p^{-2\kappa} = O(1)$ as \mbox{$n,p\rightarrow
\infty$}.
This condition imposes a large number $p$
of observation points if $\kappa$ is small. However, if the functions
$X_i$ are smooth enough such that $\kappa=1,$
then $n p^{-2\kappa} = O(1)$ is already fulfilled if $\frac{\sqrt
{n}}{p}=O(1)$ as $n,p\rightarrow\infty$,
which does not seem to be restrictive in view of practical
applications. The above result then becomes
%
\begin{equation}\label{raten}
\| \widehat{\bolds{\alpha}} - \bolds{\alpha} \|
_{\Gamma_{n,p}}^{2} =
O_{P} \bigl( \rho+ \bigl(n\rho^{{1}/({2m+2q+1})}\bigr)^{-1} \bigr).
\end{equation}
Choosing $\rho\sim n^{-(2m+2q+1)/(2m+2q+2)}$, we can conclude that
%
\begin{equation}\label{ratensmo}
\| \widehat{\bolds{\alpha}} - \bolds{\alpha} \|
_{\Gamma_{n,p}}^{2} =
O_{P} \bigl( n^{-(2m+2q+1)/(2m+2q+2)} \bigr).
\end{equation}

The next theorem studies the behavior of the estimator for the
empirical $L^{2}$-norm $\|\cdot\|_{\Gamma_{n}}$. It is shown that if $p$
is sufficiently large compared to
$n$, then based on an optimal choice of $\rho$,
the rate of convergence given in (\ref{ratensmo}) generalizes to
the semi-norm $\|\cdot\|_{\Gamma_{n}}$.

\begin{theorem} \label{thm2}
Assume \textup{(A.1)--(A.3)} as well as $n p^{-2\kappa} = O(1)$, $\rho
\rightarrow0$, $1/\break(n\rho)\rightarrow0$ as $n,p\rightarrow\infty$. Then
%
\begin{equation}\label{resthm22}
\| \widehat{\alpha} - \alpha\|_{\Gamma_n}^{2} = O_{P}
\bigl( \rho+ \bigl(n\rho^{{1}/({2m+2q+1})}\bigr)^{-1} \bigr).
\end{equation}
\end{theorem}

We finally investigate in the next theorem the behavior of $\| \widehat
{\alpha} - \alpha\|_{\Gamma}^{2}$.
The following assumption describes the additional conditions used to
derive our results.
It is well known that the covariance operator $\Gamma$ is a nuclear,
self-adjoint and nonnegative Hilbert--Schmidt operator.
We will use $\zeta_1,\zeta_2,\ldots$ to denote a complete orthonormal
system of eigenfunctions
of $\Gamma$ corresponding to the
eigenvalues $\lambda_1\ge\lambda_2\ge\cdots$.

(A.4) There exists a constant $C_7<\infty$ such that
%
\begin{eqnarray}\label{H4}
&&\operatorname{Var}\Biggl(\frac{1}{n}\sum_i
\langle X_i-\mathbb{E}(X),\zeta_r\rangle\langle X_i-\mathbb
{E}(X),\zeta_s\rangle\Biggr) \nonumber\\[-8pt]\\[-8pt]
&&\qquad\le
\frac{C_7}{n}\mathbb{E}\bigl(\langle X-\mathbb{E}(X),\zeta_r\rangle^2\bigr)
\mathbb{E}\bigl(\langle X-\mathbb{E}(X),\zeta_s\rangle^2\bigr)
\nonumber
\end{eqnarray}
holds for all $n$ and all $r,s=1,2,\ldots.$ Moreover, $\| \overline
{X}-\mathbb{E}(X))\|^2=O_P(n^{-1})$.

Relation (\ref{H4}) establishes a moment condition. It is
necessarily fulfilled if $X_1,\ldots,X_n$ are i.i.d. Gaussian random
functions. Then
$\langle X_i-\mathbb{E}(X),\zeta_r\rangle\sim
N(0,\mathbb{E}(\langle X_i-\mathbb{E}(X),\zeta_r\rangle^2))$, and
$\langle
X_i-\mathbb{E}(X),\zeta_r\rangle$ is independent of $\langle X_i-\mathbb
{E}(X),\zeta_s\rangle$
if $r\ne s$. Relation (\ref{H4}) then is an immediate consequence.

However, the validity of (\ref{H4}) does not require independence of
the functions $X_i$. For example,
in the Gaussian case, (\ref{H4}) may also be verified if
$\operatorname{Cov}( \langle X_i-\mathbb{E}(X),\zeta_r\rangle\langle
X_i-\mathbb{E}(X),\zeta_s\rangle,\langle X_j-\mathbb{E}(X),\zeta
_r\rangle\langle
X_j-\mathbb{E}(X),\zeta_s\rangle)
\le C_7\mathbb{E}(\langle
X_i-\mathbb{E}(X),\zeta_r\rangle^2)
\mathbb{E}(\langle X_i-\mathbb{E}(X),\zeta_s\rangle^2)\cdot q^{|i-j|}$
for some $0<q<1$, $C_7<\infty$ and $i\neq j$. This is of importance in our
application to ozone pollution forecasting which deals with
a \textit{time series} of functions $X_1,\ldots,X_n$.

\begin{theorem} \label{thm3}
Under the conditions of Theorem \ref{thm2} together with
assumption~\textup{(A.4)} we have
%
\begin{equation}\label{resthm3}
\| \widehat{\alpha} - \alpha\|_{\Gamma}^{2} =
O_{P} \bigl( \rho+ \bigl(n\rho^{{1}/({2m+2q+1})}\bigr)^{-1} +n^{-(2q+1)/2}
\bigr).
\end{equation}
Furthermore, (\ref{prederror}) holds for any random function $X_{n+1}$
possessing the same distribution as $X$ and independent of
$X_1,\ldots,X_n$.
\end{theorem}

Theorem \ref{thm3} shows that if $2q\ge1$ and $\rho\sim
n^{-(2m+2q+1)/(2m+2q+2)}$, then the prediction error
can be bounded by
\[
\mathbb{E} \bigl((\widehat{\alpha_0}+\langle\widehat{\alpha
},X_{n+1}\rangle-
\alpha_0-\langle\alpha,X_{n+1}\rangle
)^2| \widehat{\alpha_0},\widehat{\alpha}
\bigr)= O_{P} \bigl(n^{-({2m+2q+1})/({2m+2q+2})} \bigr).
\]

\subsection{Optimality of the rates of convergence} \label{sec32}

\noindent
For simplicity we will rely on the special case of (\ref
{funclinearregression})
with $\alpha_0=0$. In this case $\mathbb{E}
((\langle\alpha,X_{n+1}\rangle-\widehat{\alpha_0}-\langle\widehat
{\alpha},X_{n+1}\rangle)^2| \widehat{\alpha_0},\widehat{\alpha}
)\ge\| \widehat{\alpha} - \alpha\|_{\Gamma}^{2}$ if $X$
possesses a centered distribution
with $\mathbb{E}(X)=0$.
In Proposition \ref{prop1} below we then show that for suitable Sobolev spaces
of functions $\alpha$
and a large class of possible distributions of $X_i$, the rate
$n^{-(2m+2q+1)/(2m+2q+2)}$ is a \textit{lower bound}
for the rate of convergence of the prediction error over all
estimators of $\alpha$ to be computed from corresponding observations
$(X_i,Y_i)$, $i=1,\ldots,n$. Consequently,
the rate attained by our smoothing spline estimator
$\widehat{\alpha}$ must be interpreted as a minimax rate over these classes.

We first have to introduce some additional notation. For simplicity, we
will assume that the functions
$X_i(t)$ are known for all $t$ so that the number $p$ of observation
points may be chosen arbitrarily large.
We will use $\mathcal{C}_{m,D}$ to denote the space of all $m$-times
continuously differentiable functions $\alpha$ with $\int_0^1
\alpha^{(j)}(t)^2\,dt\le D$ for all $j=0,1,\ldots,m$. Furthermore, let
$\mathcal{P}_{q,C}$ denote the space of all
centered probability distributions on $L^2([0,1])$ with the properties
that (a) the sequence of eigenvalues of the
corresponding covariance operator satisfies $\sum_{j=k+1}^\infty
\lambda_j \le Ck^{-2q}$ for all sufficiently large $k$, and that (b)
the smoothing spline estimator
$\widehat{\alpha}$ satisfies $\| \widehat{\alpha} - \alpha\|_{\Gamma
}^{2} =
O_{P} (n^{-(2m+2q+1)/(2m+2q+2)})$ for $\alpha\in\mathcal{C}_{m,D}$ and
$\rho\sim n^{-(2m+2q+1)/(2m+2q+2)}$
(whenever $p$ is chosen sufficiently large compared to~$n$).
Finally, for
given $\alpha\in\mathcal{C}_{m,D}$, probability distribution $P\in
\mathcal{P}_{q,C}$ and i.i.d. random functions $X_1,\ldots,X_n$,
$X_i\sim P$, let $\hat a(\alpha,P)$ denote an arbitrary estimator of
$\alpha$ based on corresponding data $(X_i,Y_i)$, $i=1,\ldots,n$,
generated by (\ref{funclinearregression}) (with $\alpha_0=0$).

\begin{proposition}\label{prop1}
Let $c_n$ denote an arbitrary sequence of positive numbers with
$c_n\rightarrow0$ as
$n\rightarrow\infty$, and let $2q=1,3,5,\ldots.$ Under the above
assumptions, we have
\[
\lim\limits_{n\rightarrow\infty}
\sup\limits_{P\in\mathcal{P}_{q,C}} \sup\limits_{\alpha\in\mathcal{C}_{m,D}} \inf\limits_{\hat
a(\alpha,P)}\mathbb{P}\bigl( \| \alpha- \hat a(\alpha,P)\|^2_\Gamma\ge
c_n\cdot n^{-(2m+2q+1)/(2m+2q+2)}
\bigr)= 1.
\]
\end{proposition}

It is of interest to compare our
results with those of
Cai and Hall \cite{CaHa06} who analyze
the
error $\langle\alpha-\widehat{\alpha},x\rangle^2$ for a \textit{fixed} curve $x$.
Similarly to our results, the rate
of decrease of the eigenvalues $\lambda_r$ of $\Gamma$ plays an
important role. Note that, as shown
in the proof of Theorem \ref{thm3}, assumption (A.3) yields $\sum_{r=k+1}^\infty
\lambda_r=O(k^{-2q})$. Since
$\lambda_1\ge\lambda_2\ge\cdots$ this in turn implies that $\lambda
_r=O(r^{-2q-1})$, and one may
reasonably assume that $B^{-1}r^{-2q-1} \le\lambda_r \le Br^{-2q-1}$
for some $0<B<\infty$.
However,
Cai and Hall \cite{CaHa06} measure ``smoothness''
of $\alpha$ in terms of a spectral decomposition $\alpha(t)=\sum_r
\alpha_r\zeta_r(t)$
and not with respect to usual
smoothness classes. Their quantity of interest is the rate $\beta>1$
of decrease
$|\alpha_r|=O(r^{-\beta})$ as $r\rightarrow\infty$. But recall that
the error in expanding an
$m$-times continuously differentiable function with respect to $k$
suitable basis functions
(as, e.g., orthogonal polynomials or Fourier functions) is of an order of
at most $k^{-2m}$. For the sake of comparison,
assume that $\zeta_1,\zeta_2,\ldots$ define an appropriate basis for
approximating smooth functions and that
$ \inf_{f\in \operatorname{span}\{\zeta_1,\ldots,\zeta_k\}}\|\alpha-f\|^2=\sum
_{r=k+1}^\infty\alpha_r^2=O(k^{-2m})$.
This will require that $\alpha_r^2=O(r^{-2m-1})$ and, hence, $2\beta=2m+1$.

Results as derived by Cai and Hall \cite{CaHa06} additionally depend
on the spectral decomposition
$x(t)=\sum_r x_r\zeta_r(t)$ of a function $x$ of interest.
The essential condition on the structure of the coefficients $x_r$ may
be re-expressed in the
following form: There exist some $\nu\in\mathbb{R}$ and $0<D_0<\infty$
such that
$D_0^{-1}r^{\nu} \le\frac{x_r^2}{\lambda_r} \le D_0r^{\nu}$ for all
$r=1,2,\ldots.$
Rates of convergence then follow from the magnitude of
$\nu$, and it is shown that parametric rates $n^{-1}$ (or $n^{-1}\log
n $) are achieved
if $\nu\le-1$.

Now consider a random function $X_{n +1}$ and assume that the
underlying distribution is Gaussian. It is then
well known that
$X_{n+1}(t)=\sum_r x_{n+1,r} \zeta_r(t)$ for
independent $N(0,\lambda_r)$-distributed coefficients $x_{n+1,r}$.
Consequently, $\frac{x_{n+1,r}^2}{\lambda_r}$
are i.i.d. $\chi_1^2$-distributed variables for all $r=1,2,\ldots,$ and
if $\nu\le0$ we obtain
$ \mathbb{P}(D_0^{-1}r^{\nu} \le\frac{x_{n+1,r}^2}{\lambda_r}
\le D_0r^{\nu}\mbox{ for all } r=1,2,\ldots)=0$
for all $0<D_0<\infty$. This already shows that parametric rates
$n^{-1}$ cannot be achieved for the error $\langle\alpha-\widehat
{\alpha},X_{n+1}\rangle^2$. On the other hand, for arbitrary $\nu>0$
and $0<\delta<1$ we have
$ \mathbb{P}(D_0^{-1}r^{\nu} \le\frac{x_{n+1,r}^2}{\lambda_r}
\le D_0r^{\nu}\mbox{ for all } r=1,2,\ldots)
\ge\delta$, whenever $D_0$ is sufficiently large. If $B^{-1}r^{-2q-1}
\le\lambda_r \le Br^{-2q-1}$ and
$\alpha_r^2=O_P(r^{-2m+1})$,
then
for a function $x$ with
$D_0^{-1}r^{\nu} \le\frac{x_{n+1,r}^2}{\lambda_r} \le D_0r^{\nu}$,
$\nu>0$,
the convergence rates of Cai and Hall \cite{CaHa06}
translate into
\[
\langle\widehat{\alpha}-\alpha,x\rangle^2 =O_P\bigl(n^{-(2m +2q+1-2\nu)/(2m+2q+2)}\bigr),
\]
which provides an additional motivation for the fact that
the rates derived in our paper constitute a lower bound.
For non-Gaussian distributions a comparison is more difficult, since
under assumption (A.4) only the Chebyshev inequality
may be used to bound the probabilities $D_0^{-1}r^{\nu} \le\frac
{x_{n+1,r}^2}{\lambda_r} \le D_0r^{\nu}$.

Another statistically very different problem consists in an optimal
estimation of~$\alpha$ by $\widehat{\alpha}$ with
respect to the usual $L^2$-norm. In a recent work, Hall and
Horowitz~\cite{HaHo07} derive optimal rates of convergence
of $\| \widehat{\alpha}-\alpha\|^2$. These rates again depend on the
rate of decrease
$|\alpha_r|=O(r^{-\beta})$.
Recall that our assumptions do not provide any link between $\alpha$
and $X_i$;
part of the structure of $\alpha$ may not even be identifiable.
Indeed, under assumptions (A.1)--(A.4) there is no way to guarantee
that the bias
$\| \alpha-\mathbb{E}_\varepsilon(\widehat{\alpha})\|^2$
converges to zero and it can only be shown that
$\| \widehat{\alpha}-\alpha\|^2=O_P(1)$ (see the proof of Theorem
\ref{thm2} below). This already highlights the
theoretical difference between optimal estimation with respect to $\|
\widehat{\alpha}-\alpha\|_\Gamma^2$
and $\| \widehat{\alpha}-\alpha\|^2$. Based on additional assumptions
as indicated above,
although sensible bounds for the bias may be derived, it must be
emphasized that an estimator minimizing
$\| \widehat{\alpha}-\alpha\|^2$ will have to rely on $\rho\gg
n^{-(2m+2q+1)/(2m+2q+2)}$, which corresponds
to an oversmoothing with respect to
$\| \widehat{\alpha}-\alpha\|_\Gamma^2$.
This effect has already been noted by Cai and Hall \cite{CaHa06}. In
our context,
without additional assumptions linking the eigenvalues of $\Gamma$ and
of the spline matrix $\mathbf{A}_{m}$,
the only general bound
for the $L_2$-variability of the estimator is
$\| \widehat{\alpha}-\mathbb{E}_\varepsilon(\widehat{\alpha})\|
^2=O_P(\frac{1}{n\rho})$ (this
result may be derived by arguments similar to those used in the proofs
of our theorems). With
$\rho= n^{-(2m+2q+1)/(2m+2q+2)}$ this leads to
$\| \widehat{\alpha}-\mathbb{E}_\varepsilon(\widehat{\alpha})\|
^2=O_P(n^{-1/(2m+2q+2)})$, and
better rates may only be achieved with $\rho\gg
n^{-(2m+2q+1)/(2m+2q+2)}$. A more detailed study of
this problem is not in the scope of the present paper.

\subsection{Choice of smoothing parameters} \label{sec33}

The above result of Section \ref{sec31} implies that the
choice of the smoothing parameter $\rho$ is of crucial importance.
A natural way to determine $\rho$ is to minimize a leave-one-out
cross-validation criterion. We preferably adapt the simplified
\textit{Generalized Cross-Validation} (GCV) introduced by Wahba \cite
{Wa90} in the context
of smoothing splines. For fixed $m$, in our application the
GCV criterion takes the form
%
\begin{equation}\label{gcvcrit}
\mathit{GCV}_m(\rho)
:=\frac{({1}/{n})\|\mathbf{Y}-\mathbf{H}_\rho
\mathbf{Y}\|^2}{(1-n^{-1}\operatorname{Tr}(\mathbf{H}_\rho))^2},
\end{equation}
where
$\mathbf{H}_\rho:=(np)^{-1}\mathbf{X} (\frac{1}{n
p^2}\mathbf{X}^\tau\mathbf{X}+\frac{\rho}{p}\mathbf{A}_m)^{-1}
\mathbf{X}^\tau$.

Proposition \ref{prop2} below provides a justification for the use of the GCV
criterion. Recall that the estimators
$\widehat{\bolds{\alpha}}\equiv
\widehat{\bolds{\alpha}}_{\rho;m}$ depend on $\rho$ as well as
on the spline order~$m$. Obviously,
$\frac{1}{p}\mathbf{X}\widehat{\bolds{\alpha}}_{\rho;m}=\mathbf
{H}_\rho
\mathbf{Y}$ is an estimator
of the conditional mean $(\langle X_1-\overline{X},\alpha\rangle,\ldots
,\langle X_n-\overline{X},\alpha\rangle
)^\tau$ of $\mathbf{Y}$ given $X_1,\ldots,X_n$. Let
\[
\mathit{ASE}_m(\rho):=\frac{1}{n}\sum_i\Biggl[\langle X_i-\overline{X},\alpha
\rangle-
\frac{1}{p}\sum_j \bigl(X_i(t_j)-\overline{X}(t_j)\bigr)\hat
\alpha_{\rho;m}(t_j)\Biggr]^2
\]
denote the average squared error of this estimator. The only difference between
$\mathit{ASE}_m(\rho)$ and $\|\widehat{\bolds{\alpha}}_{ \rho}-
\bolds{\alpha} \|_{\Gamma_{n,p}}^{2}$ is the discretization error
encountered when approximating
$\langle X_i,\alpha\rangle$ by $\frac{1}{p}\sum_j X_i(t_j) \alpha
(t_j)$, and hence
$\mathit{ASE}_m(\rho)=\|\widehat{\bolds{\alpha}}_{ \rho}-
\bolds{\alpha} \|_{\Gamma_{n,p}}^{2}+O_P(p^{-2\kappa})$.

If $\hat\rho$ denotes the minimizer of GCV for fixed $m$, we can
conclude from relation (\ref{gcvrho}) of Proposition \ref{prop2} that the
error $\mathit{ASE}_m(\hat\rho)$ is asymptotically first-order equivalent to
the error
$\mathit{ASE}_m(\rho_{\mathrm{opt}})$ to be obtained from an optimal choice of the
smoothing parameter.
Furthermore, (\ref{gcvrhom}) shows that an analogous result holds
if GCV is additionally used to select the order $m$ of the smoothing
spline, which means that the optimal rate can be reached adaptively.

\begin{proposition}\label{prop2}
In addition to assumptions \textup{(A.1)--(A.3)} as well as $n
p^{-2\kappa}=O(1),$ suppose that
$\mathbb{E}(\exp(\beta\varepsilon_i^2))<\infty$ for some $\beta>0$.
If for fixed $m$, $\hat\rho$ denotes the minimizer of $\mathit{GCV}(\rho)$
over $\rho\in[n^{-2m+\delta},\infty)$ for some $\delta>0$, then
%
\begin{equation}\label{gcvrho}
|\mathit{ASE}_m(\hat\rho)-\mathit{ASE}_m(\rho_{\mathrm{opt}})|=O_P\bigl(n^{-{1}/{2}}
\mathit{ASE}_m(\rho_{\mathrm{opt}})^{{1}/{2}}\bigr),
\end{equation}
where $\rho_{\mathrm{opt}}$ minimizes
$\mathit{MSE}_m(\rho):=\mathbb{E}_\varepsilon(\mathit{ASE}_m(\rho))$ over all $\rho>0$.

Furthermore, if $\hat m,\hat\rho$ denotes the minimizers of (\ref
{gcvcrit}) over $\rho\in[n^{-2m+\delta},\infty)$, $\delta>0$, and
$m=1,\ldots,M_n$, $M_n\le n/2$, then
%
\begin{equation}\label{gcvrhom}\qquad
|\mathit{ASE}_{\hat
m}(\hat\rho)-\mathit{ASE}_{m_{\mathrm{opt}}}(\rho_{\mathrm{opt}})|=O_P(n^{-{1}/{2}}
\mathit{ASE}_{m_{\mathrm{opt}}}(\rho_{\mathrm{opt}})^{{1}/{2}}\log M_n),
\end{equation}
where $\rho_{\mathrm{opt}},m_{\mathrm{opt}}$ minimize
$\mathit{MSE}_m(\rho):=\mathbb{E}_\varepsilon(\mathit{ASE}_m(\rho))$ over all $\rho>0$ and
$m=1,\ldots,M_n$.
\end{proposition}

\section{Case of a noisy covariate} \label{sec4}

In a number of important applications measurements of the
explanatory curves $X_i$ may be contaminated by noise. There then
additionally exists
an errors-in-variable problem complicating further analysis.
Our setup is inspired by other works dealing with noisy observations of
functional
data (e.g., Cardot \cite{Ca00} or Chiou, M\"uller and Wang \cite{ChMuWa03}):
At each point $t_j$ the corresponding functional value
$X_i(t_j)$ is corrupted by some random error $\delta_{ij}$ so that
actual observations
$W_i(t_j)$ are given by
%
\begin{equation}\label{noiseXdiscr}
W_{i} (t_{j}) = X_{i} (t_{j}) + \delta_{ij},\qquad i=1,\ldots,n, j=1,\ldots
,p,
\end{equation}
where $(\delta_{ij})_{i=1,\ldots,n,j=1, \ldots, p}$ is a
sequence of independent real random variables such that
for all $i=1,\ldots,n$ and all $j=1, \ldots, p$
%
\begin{equation}\label{conddelta}
\mathbb{E}_{\varepsilon}(\delta_{ij}) = 0,\qquad
\mathbb{E}_{\varepsilon} (\delta_{ij}^{2}) = \sigma_{\delta}^{2}
\quad\mbox{and}\quad\mathbb{E}_\varepsilon(\delta_{ij}^4) \le C_8
\end{equation}
for some constant $C_8>0$ (independent of $n$ and $p$).
We furthermore assume that~$\delta_{ij}$ is independent of $\varepsilon_i$
and of the $X_i$'s.

In this situation,
an analogue of our estimator $\widehat\alpha$ of Section \ref{sec2} can
still be computed by replacing in (\ref{solX})
the (unknown) matrix $\mathbf{X}$ by the $n \times p$ matrix
$\mathbf{W}$ with general
terms $W_{i} (t_{j})-\overline{W}$, $i=1,\ldots,n$, $j=1,\ldots,p$.
However, performance of the resulting estimator will suffer from the
additional noise in the observations. If the error variance $\sigma
^2_\delta$ is large, there may exist a substantial difference
between $\mathbf{X}^{\tau}\mathbf{X}$
and $\mathbf{W}^{\tau}\mathbf{W}$. Indeed, $\mathbf
{W}^{\tau}\mathbf{W}$ is a biased estimator
of $\mathbf{X}^{\tau}\mathbf{X}$:
%
\begin{equation}\label{debruitage}
\frac{1}{n p^{2}} \mathbf{W}^{\tau} \mathbf{W} = \frac{1}{n p^{2}}
\mathbf{X}^{\tau} \mathbf{X} +
\frac{\sigma_{\delta}^{2}}{p^{2}} \mathbf{I}_{p} + \mathbf{R},
\end{equation}
where $\mathbf{R}$ is a $p\times p$ matrix such that its
largest singular value is of
order $O_{P} ( \frac{1}{n^{1/2} p} )$, (see the proof of
Theorem \ref{thm4} below).
This result suggests that we use $\frac{1}{n p^{2}} \mathbf{W}^{\tau}
\mathbf{W} - \frac{\sigma_{\delta}^{2}}{p^{2}} \mathbf{I}_{p}$ as an
approximation of $\frac{1}{n p^{2}} \mathbf{X}^{\tau} \mathbf{X}$. A
prerequisite is, of course, the availability of
an estimator $\hat\sigma_{\delta}^{2}$
of the unknown variance $\sigma_{\delta}^{2}$. Following
Gasser, Sroka and Jennen-Steinmetz \cite{GaSrSt86}, we will rely on
%
\begin{equation}\label{sigma2hat}\quad
\widehat{\sigma}_{\delta}^{2} := \frac{1}{n} \sum_{i=1}^{n} \frac
{1}{6(p-2)} \sum_{j=2}^{p-1}
[ W_{i} (t_{j-1}) - W_{i} (t_{j}) + W_{i} (t_{j+1}) - W_{i}
(t_{j}) ]^{2}.
\end{equation}
These arguments now lead to the following modified estimator
$\widehat{\bolds{\alpha}}_{\mathbf{W}}$ of
$\bolds{\alpha}$ in the case of noisy observations:
%
\begin{equation}\label{solssftlsbis}
\widehat{\bolds{\alpha}}_{\mathbf{W}} := \frac{1}{n p} \biggl(
\frac{1}{n p^{2}}
\mathbf{W}^{\tau} \mathbf{W} + \frac{\rho}{p} \mathbf{A}_{m}
- \frac{\widehat{\sigma}_{\delta}^{2}}{p^{2}} \mathbf{I}_{p} \biggr)^{-1}
\mathbf{W}^{\tau} \mathbf{Y}.
\end{equation}
An estimator of the function $\alpha$ is given
by $\widehat{\alpha}_{\mathbf{W}} = s_{\widehat{\bolds{\alpha
}}_{\mathbf{W}}}$,
where $s_{\widehat{\bolds{\alpha}}_{\mathbf{W}}}$ is again
the natural spline interpolant of order $2m$
as defined in Section \ref{sec2}.

We want to note that $\widehat{\bolds{\alpha}}_{\mathbf{W}}$
is closely related to an estimator proposed by
Cardot et al. \cite{CaCrKnSa07}. The latter is motivated by
the \textit{Total Least Squares}
(TLS) method (see, e.g., Golub and Van Loan \cite{GoVanLo80},
Fuller \cite{Fu87}, or Van Huffel and Vandewalle \cite{VanHuVa91}) and
the only difference from (\ref{solssftlsbis}) consists in the use
of a correction term slightly different from $- \frac{\widehat{\sigma
}_{\delta}^{2}}{p^{2}} \mathbf{I}_{p}$.

Of course there are many alternative strategies for dealing with the
errors-in-variable problem induced by (\ref{noiseXdiscr}).
A straightforward approach, which is \mbox{frequently} used in functional
data analysis, is to apply nonparametric smoothing procedures
in order to obtain estimates $\hat X_i(t_j)$ from the data
$(W_i(t_j),t_j)$. When replacing
$\mathbf{X}$ by $\widehat{\mathbf{X}}$ in (\ref{solX}),
one can then define a ``smoothed'' estimator $\widehat{\bolds
{\alpha}}_{S}$.
Of course this estimator may be as efficient as (\ref{solssftlsbis}),
but it is computationally more involved and
appropriate smoothing parameters have to be selected for nonparametric
estimation of each curve $X_i$.

Our aim is now to study the asymptotic behavior of $\widehat{\alpha
}_{\mathbf{W}}$. Theorem \ref{thm4} below provides bounds
(with respect to the semi-norm $\Gamma_{n,p}$) for the difference
between~$\widehat{\bolds{\alpha}}_{\mathbf{W}}$
and the ``ideal'' estimator $\widehat{\bolds{\alpha}}$ defined for
the true curves $X_1,\ldots,X_n$.
We will impose the following additional condition on the function
$\alpha$:
\begin{enumerate}[(A.5)]
\item[(A.5)] For every $\delta>0$ there exists a constant $C_\alpha
<\infty$ such that
\[
\frac{1}{p^{1/2}}\biggl\|\frac{1}{n p} \mathbf{X}^\tau\mathbf
{X}\bolds{\alpha}\biggr\| > C_\alpha,
\]
holds with probability larger or equal to $1-\delta$.
\end{enumerate}

\begin{theorem} \label{thm4}
Assume \textup{(A.1), (A.2), (A.5)} as well as $n p^{-2\kappa}=O(1)$, $\rho
\rightarrow0$,
$1/(n \rho) \rightarrow0$ as $n,p \rightarrow\infty$. Then
%
\begin{equation}\label{resthm4}
\| \widehat{\bolds{\alpha}}_{\mathbf{W}} - \widehat
{\bolds{\alpha}} \|_{\Gamma_{n,p}}^2 = O_{P} \biggl( \frac
{1}{n p \rho} + \frac{1} {n}\biggr).
\end{equation}
\end{theorem}

Together with assumption (A.3) we can therefore conclude from Theorems
\ref{thm1} and \ref{thm4} that
\[
\| \widehat{\bolds{\alpha}}_{\mathbf{W}} - \bolds{\alpha
} \|_{\Gamma_{n,p}}^{2} =
O_{P} \biggl( \rho+ \bigl( n \rho^{{1}/({2m+2q+1})}\bigr)^{-1} +
\frac{1}{n p \rho} \biggr).
\]

We have already seen in Section \ref{sec3} that the optimal order of
the two first terms is reached for a choice
of $\rho\sim n^{-(2m+2q+1)/(2m+2q+2)}$.
From an asymptotic point of view, the use of $\widehat{\alpha}_{\mathbf
{W}}$ results in the addition of the extra
term $1/(n p \rho)$ in the rate of convergence. For $\rho\sim
n^{-(2m+2q+1)/(2m+2q+2)}$ we
have $1/(n p \rho) \sim n^{-1/(2m+2q+2)}/p$. This term is of order
$n^{-(2m+2q+1)/(2m+2q+2)}$ for
$p \sim n^{(2m+2q-1)/(2m+2q+2)}$. This means that the $\widehat{\alpha
}_{\mathbf{W}}$ reaches the same
rate of convergence as $\widehat{\bolds{\alpha}}$ provided that
$p$ is sufficiently large
compared to $n$. More precisely, it is required that $p\geq C_p \max
(n^{1/2\kappa}, n^{(2m+2q-1)/(2m+2q+2)})$ for some positive constant $C_p$.

As shown in Theorem \ref{thm5} below, these qualitative results generalize when
considering the semi-norms $\Gamma_{n}$ or $\Gamma$.

\begin{theorem} \label{thm5}
Assume \textup{(A.1)--(A.3), (A.5)} as well as $np^{-2\kappa} =O(1)$, $\rho
\rightarrow0$, $1/(n \rho) \rightarrow0$ as $n,p \rightarrow\infty$. Then
%
\begin{equation}\label{resthm5}
\| \widehat{\alpha}_{\mathbf{W}} - \widehat{\alpha} \|
_{\Gamma_n}^{2} = O_{P} \biggl( \frac{1}{n p \rho} + \frac{1}{n} \biggr),
\end{equation}
and if assumption \textup{(A.4)} is additionally satisfied,
%
\begin{equation}\label{resthm6}
\| \widehat{\alpha}_{\mathbf{W}} - \widehat{\alpha} \|
_{\Gamma}^{2} =
O_{P} \biggl( \frac{1}{n p \rho} + \frac{1}{n} + n^{-(2q+1)/2} \biggr).
\end{equation}
\end{theorem}

\section{Application to ozone pollution forecasting} \label{sec5}

In this section, our methodology is applied to the problem of
predicting the level of ozone pollution.
For our analysis, we use a data set collected by
ORAMIP (Observatoire Régional de l'Air en Midi-Pyrénées){}, an air
observatory located in the
city of Toulouse (France). The
concentration of specific pollutants as well as meteorological
variables are measured each hour. Some
previous studies using the same data are described in
Cardot, Crambes and Sarda \cite{CaCrSa06} and
Aneiros-Perez et al. \cite{AnCaEsVi04}.

The response variable $Y_i$ of interest is the maximum of ozone for a
day. Repeated measurements
of ozone concentration obtained for the \textit{preceding} day are used
as a functional explicative variable
$X_i$. More precisely, each $X_i$ is observed at $p=24$ equidistant
points corresponding to hourly measurements.
The sample size is $n = 474$. It is assumed that the relation between
$Y_i$ and $X_i$
can be modeled by the functional linear regression model (\ref
{funclinearregression}).
We note at this point that $X_1,X_2,\ldots$ constitute a time series of
functions, and that it is therefore reasonable
to suppose some correlation between the $X_i$'s. The results of an
earlier, unpublished study indicate that there only
exists some
``short memory'' dependence.

Now, for a curve $X_{n+1}$ outside the sample, we want to predict
$Y_{n+1}$, the maximum of ozone the day after.
Assuming that $(X_{n+1},Y_{n+1})$ follows the same model~(\ref
{funclinearregression}) and using our estimators $\widehat\alpha$
of $\alpha$ and $\widehat{\alpha_0}$ of $\alpha_0$ described in
Section \ref{sec2}, a~predictor
$\widehat{Y}_{n+1}$ is given by the formula
%
\begin{equation}\label{hatY}
\widehat{Y}_{n+1} := \widehat{\alpha_0}+\int_I\widehat{\alpha} (t)
X_{n+1} (t) \,dt.
\end{equation}
It cannot be excluded that actual observations of $X_i$ may
be contaminated with noise. We will thus additionally
consider the modified estimator $\widehat{\alpha}_{\mathbf{W}}$
developed in Section \ref{sec4} and the corresponding predictor
$\widehat{Y}_{\mathbf{W},n+1}$. For simplicity, the integral in~(\ref{hatY})
is approximated by $\frac{1}{p} \sum_{j=1}^{p} \widehat{\alpha}
(t_{j}) X_{n+1} (t_{j})$.
With additional assumptions on the $\varepsilon_i$'s we can also build
asymptotic intervals
of prediction for $Y_{n+1}$. Indeed, let us assume that $\varepsilon
_1,\ldots,\varepsilon_{n+1}$ are i.i.d.
random variables having a normal distribution $\mathcal{N} ( 0 ,
\sigma_{\varepsilon}^{2} )$. The first point is
to estimate the residual variance $\sigma_{\varepsilon}^{2}$. A
straightforward estimator is given by the empirical variance
%
\begin{equation}\label{hatsigmaeps2}
\widehat{\sigma_{\varepsilon}}^{2} := \frac{1}{n} \sum_{i=1}^{n}
\Biggl( Y_{i} - \overline{Y}-\frac{1}{p} \sum_{j=1}^{p} \widehat{\alpha
} (t_{j}) \bigl(X_{i} (t_{j})-\overline{X}(t_j)\bigr) \Biggr)^{2}.
\end{equation}
Our theoretical results imply that $\widehat{\sigma_{\varepsilon}}$ is a
consistent estimator of $\sigma_{\varepsilon}^{2}$. Furthermore,
we can then infer from Theorem \ref{thm3} that
$\frac{Y_{n+1}-\widehat{Y}_{n+1}}{\widehat{\sigma_{\varepsilon}}}$
asymptotically follows a standard normal distribution.
Given $\tau\in\,]
0,1[$, an asymptotic $(1 - \tau)$-prediction interval for~$Y_{n+1}$ can
be derived as
%
\begin{equation}\label{inthatY}
[ \widehat{Y}_{n+1} - z_{1 - \tau/ 2} \widehat{\sigma_{\varepsilon}}
, \widehat{Y}_{n+1} + z_{1 - \tau/ 2} \widehat{\sigma_{\varepsilon}}
],
\end{equation}
where $z_{1 - \tau/ 2}$ is the quantile of order $1 - \tau/
2$ of the $\mathcal{N} (0,1)$ distribution.
Of course, the same developments are valid when one replaces $\widehat
{Y}_{n+1}$ by $\widehat{Y}_{\mathbf{W},n+1}$.

In order to study performance of our estimators we
split the initial sample into two sub-samples:

%
\begin{figure}[b]

\includegraphics{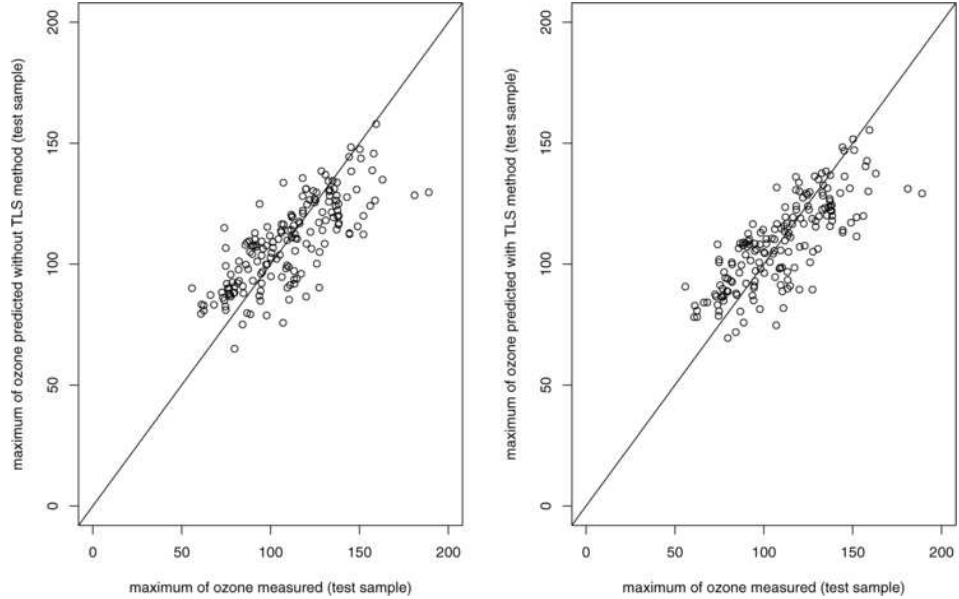}

\caption{Daily predicted values $\widehat{Y}$ (left) and
$\widehat{Y}_{\mathbf{W}}$ (right) of the maximum of ozone versus the
measured values.}\label{fig1}
\end{figure}

\begin{figure}

\includegraphics{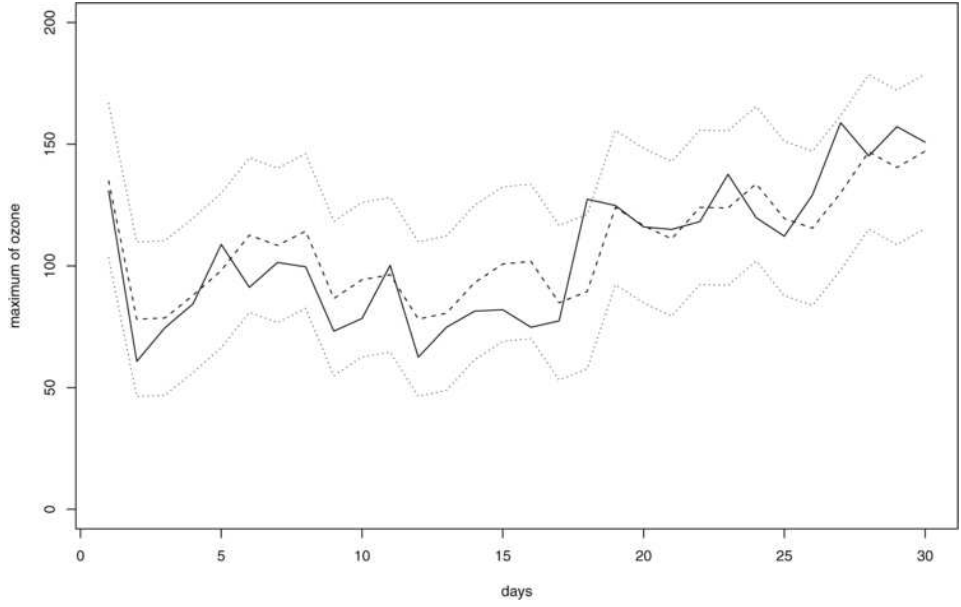}

\caption{Measured values of the maximum of ozone (solid line),
predicted values (dashed line)
and $95\%$ prediction band (dotted lines).}\label{fig2}
\end{figure}

\begin{itemize}
\item A learning sample, $(X_{i},Y_{i})_{i=1, \ldots, n_{l}}$, $n_{l} =
300$, was
used to determine the estimators $\widehat{\alpha}$ and $\widehat
{\alpha}_{\mathbf{W}}$.
\item A test sample, $(X_{i},Y_{i})_{i=n_{l} + 1, \ldots, n_{l} +
n_{t}}$, $n_{t} = 174$, was used
to evaluate the quality of the estimation.
\end{itemize}
Construction of estimators was based on $m=2$ (cubic smoothing
splines), and
the smoothing parameters $\rho$ were selected by minimizing $\mathit{GCV}(\rho)$
as defined in~(\ref{gcvcrit}). Note\vspace*{-5pt} that GCV for $\widehat{\alpha}_{\mathbf{W}}$
requires that the matrix
$\frac{1}{n p^2}\mathbf{X}^\tau\mathbf{X}$ in the definition of $\mathbf
{H}_\rho$
has to be replaced by
$\frac{1}{n p^{2}}
\mathbf{W}^{\tau} \mathbf{W}
- \frac{\widehat{\sigma}_{\delta}^{2}}{p^{2}} \mathbf{I}_{p}$.
Figure \ref{fig1} presents
the daily predicted values $\widehat{Y}$ and $\widehat{Y}_{\mathbf
{W}}$ of the maximum of ozone versus the measured
$Y$-values of the test sample.
Both graphics are close, which is confirmed by the computation of the
prediction error given by
\[
\mathit{EQM} ( \widehat{\alpha} ) := \frac{1}{n_{t}} \sum_{i=n_{l} +
1}^{n_{l} + n_{t}} ( Y_{i} - \widehat{Y}_{i} )^{2},
\]
with a similar definition for $\widehat{\alpha}_{\mathbf
{W}}$. We have, respectively,
$\mathit{EQM} ( \widehat{\alpha} ) = 281.97$ and $\mathit{EQM} ( \widehat
{\alpha}_{\mathbf{W}} ) = 270.13,$
which shows a very minor
advantage of the estimator $\widehat{\alpha}_{\mathbf{W}}$.
In any case, in Figure \ref{fig1} the points seem to be reasonably spread around
the diagonal $\hat Y=Y$, and the plots
do not indicate any major problem with our estimators. Corresponding
prediction intervals are given in Figure \ref{fig2}.

%
%

\section{Proof of the results} \label{sec6}

\subsection[Proof of Theorem 1]{Proof of Theorem \protect\ref{thm1}}

First consider relation (\ref{resthm1}), and note that
\[
\mathbb{E}_{\varepsilon} ( \widehat{\bolds{\alpha}} ) =
\frac{1}{n p^{2}} \biggl( \frac{1}{n p^{2}} \mathbf{X}^{\tau} \mathbf{X}
+ \frac{\rho}{p} \mathbf{A}_{m} \biggr)^{-1} \mathbf{X}^{\tau} \mathbf
{X} \bolds{\alpha}+ \frac{1}{n p} \biggl( \frac{1}{n p^{2}} \mathbf
{X}^{\tau} \mathbf{X} + \frac{\rho}{p} \mathbf{A}_{m} \biggr)^{-1}
\mathbf{X}^{\tau} \mathbf{d},
\]
where $\mathbf{d}=(d_1-\overline{d},\ldots,d_n-\overline{d})^\tau$.

It follows that $\mathbb{E}_{\varepsilon} ( \widehat
{\bolds{\alpha}} )$ is a solution of the minimization
problem
\[
\min_{\mathbf{a} \in\mathbb{R}^{p}} \biggl\{ \frac{1}{n} \biggl\|
\frac{1}{p} \mathbf{X} \bolds{\alpha} +\mathbf{d}- \frac{1}{p}
\mathbf{X}
\mathbf{a} \biggr\|^{2} + \frac{\rho}{p} \mathbf{a}^{\tau}
\mathbf{A}_{m} \mathbf{a} \biggr\}.
\]
This implies
\[
\frac{1}{n} \biggl\| \frac{1}{p} \mathbf{X} \bolds{\alpha
}+\mathbf{d} -
\frac{1}{p} \mathbf{X} \mathbb{E}_{\varepsilon} ( \widehat
{\bolds{\alpha}} ) \biggr\|^{2} +
\frac{\rho}{p} \mathbb{E}_{\varepsilon} ( \widehat{\bolds{\alpha
}} )^{\tau}
\mathbf{A}_{m} \mathbb{E}_{\varepsilon} ( \widehat{\bolds{\alpha
}} ) \leq\frac{\rho}{p} \bolds{\alpha}^{\tau}
\mathbf{A}_{m} \bolds{\alpha}+\frac{1}{n}\|\mathbf
{d}\|^{2}.
\]
But definition of $\mathbf{A}_{m}$ and (\ref{natsmooth})
lead to
\[
\frac{1}{p} \bolds{\alpha}^{\tau} \mathbf{A}_{m} \bolds{\alpha
} =
\frac{1}{p} \bolds{\alpha}^{\tau} \mathbf{P}_{m} \bolds
{\alpha} +
\int_{0}^1 s_{\bolds{\alpha}}^{(m)} (t)^{2} \,dt \leq\frac{1}{p}
\bolds{\alpha}^{\tau} \mathbf{P}_{m}
\bolds{\alpha} + \int_{0}^1 \alpha^{(m)} (t)^{2} \,dt
\]
and (\ref{resthm1}) is an immediate consequence.
Let us now consider relation (\ref{resthm12}). There exists a
complete orthonormal system of eigenvectors $u_1,u_2,\ldots,u_p$
of $ \frac{1}{n p} \mathbf{X}^{\tau} \mathbf{X}$ such that
$\frac{1}{n p} \mathbf{X}^{\tau} \mathbf{X}=\sum_{j=1}^p \lambda_{x,j}
u_ju_j^{\tau}$. Let
$k:=[\rho^{-{1}/({2m+2q+1})}]$.
By our assumptions we obtain
\begin{eqnarray}\label{thm1step1}
&& \mathbb{E}_{\varepsilon} \bigl( \|\widehat{\bolds{\alpha}} -
\mathbb{E}_{\varepsilon} ( \widehat{\bolds{\alpha}}) \|_{\Gamma
_{n,p}}^{2} \bigr) \nonumber\\
&&\qquad =
\frac{1}{p} \mathbb{E}_{\varepsilon} \biggl( \frac{1}{n^{2} p^{2}}
\bolds{\varepsilon}^{\tau}
\mathbf{X}\biggl( \frac{1}{n p^{2}} \mathbf{X}^{\tau} \mathbf{X} + \frac
{\rho}{p} \mathbf{A}_{m} \biggr)^{-1} \nonumber\\
&&\hspace*{59.1pt}{}\times \frac{1}{n p} \mathbf{X}^{\tau}
\mathbf{X} \biggl( \frac{1}{n p^{2}} \mathbf{X}^{\tau} \mathbf{X} + \frac
{\rho}{p} \mathbf{A}_{m} \biggr)^{-1} \mathbf{X}^{\tau} \bolds
{\varepsilon} \biggr) \nonumber\\
&&\qquad \leq \frac{\sigma_{\varepsilon}^{2}}{n} \operatorname{Tr} \biggl[ \biggl( \frac
{1}{n p} \mathbf{X}^{\tau} \mathbf{X}
+ \rho\mathbf{A}_{m} \biggr)^{-1} \frac{1}{n p} \mathbf{X}^{\tau}
\mathbf{X} \biggr] \\
&&\qquad = \frac{\sigma_{\varepsilon}^{2}}{n} \operatorname{Tr} \biggl[ \biggl( (\rho
\mathbf{A}_{m})^{-{1}/{2}}\biggl(\frac{1}{n p}
\mathbf{X}^{\tau} \mathbf{X}\biggr)(\rho\mathbf{A}_{m})^{-{1}/{2}} +
\mathbf{I}_{p} \biggr)^{-1} \nonumber\\
&&\hspace*{99pt}{}\times (\rho\mathbf{A}_{m})^{-{1}/{2}}\biggl(\frac
{1}{n p} \mathbf{X}^{\tau} \mathbf{X}\biggr)
(\rho\mathbf{A}_{m})^{-{1}/{2}} \biggr] \nonumber\\
&&\qquad \leq \frac{\sigma_{\varepsilon}^{2}}{n} \operatorname{Tr} ( \mathbf
{D}_{1,\rho}+
\mathbf{D}_{2,\rho} ),\nonumber
\end{eqnarray}
where
\begin{eqnarray*}
\mathbf{D}_{1,\rho}
:\!&=&\Biggl( (\rho\mathbf{A}_{m})^{-{1}/{2}}
\Biggl(\sum_{j=1}^k\lambda_{x,j} u_ju_j^{\tau}\Biggr)(\rho\mathbf{A}_{m})^{-
{1}/{2}} + \mathbf{I}_{p} \Biggr)^{-1}\\
&&{}\times
(\rho\mathbf{A}_{m})^{-{1}/{2}}\Biggl(\sum_{j=1}^k\lambda_{x,j}
u_ju_j^{\tau}\Biggr)(\rho\mathbf{A}_{m})^{-{1}/{2}}
\end{eqnarray*}
and
\begin{eqnarray*}
\mathbf{D}_{2,\rho}
:\!&=&\Biggl( (\rho\mathbf{A}_{m})^{-{1}/{2}}
\Biggl(\sum_{j=k+1}^p\lambda_{x,j} u_ju_j^{\tau}\Biggr)(\rho\mathbf
{A}_{m})^{-{1}/{2}} + \mathbf{I}_{p} \Biggr)^{-1}\\
&&{}\times
(\rho\mathbf{A}_{m})^{-{1}/{2}}\Biggl(\sum_{j=k+1}^p\lambda_{x,j}
u_ju_j^{\tau}\Biggr)(\rho\mathbf{A}_{m})^{-{1}/{2}}
\end{eqnarray*}
which are symmetric $p\times p$ matrices with
%
\begin{equation}\label{thm1step2}
\sup_{\Vert\mathbf{v} \Vert= 1} \mathbf{v}^{\tau} \mathbf
{D}_{1,\rho} \mathbf{v} < 1
\quad\mbox{and}\quad
\sup_{\Vert\mathbf{v} \Vert= 1} \mathbf{v}^{\tau} \mathbf
{D}_{2,\rho} \mathbf{v} < 1.
\end{equation}
Furthermore, $\mathbf{D}_{1,\rho}$ is of rank $k$ and therefore only
possesses $k$ nonzero eigenvalues.
Hence
%
\begin{equation}\label{thm1step3}
\operatorname{Tr} ( \mathbf{D}_{1,\rho})\le k.
\end{equation}
Let $\mathbf{a}_{1,p},\ldots,\mathbf
{a}_{m,p},\mathbf{a}_{m+1,p},\ldots,\mathbf{a}_{p,p}$ denote a
complete, orthonormal system of eigenvectors of $\mathbf{A}_m$
corresponding to the eigenvalues $\mu_{1,p}=\cdots=\mu_{m,p}=1$ and $\mu
_{m+1,p}\le\cdots\le\mu_{p,p}$.
By (\ref{thm1step1}), (\ref{thm1step2}) and (\ref{thm1step3}) as well
as (\ref{egams}), we thus obtain
\begin{eqnarray}\label{thm1step4}\qquad\quad
&&\mathbb{E}_{\varepsilon} \bigl( \|\widehat{\bolds{\alpha}} -
\mathbb{E}_{\varepsilon}
( \widehat{\bolds{\alpha}})
\|_{\Gamma_{n,p}}^{2} \bigr)\nonumber\\
&&\qquad\le\frac{\sigma_{\varepsilon}^{2}}{n}\Biggl( k+\sum_{j=1}^p \mathbf
{a}_{j,p}^{\tau} \mathbf{D}_{2,\rho}\mathbf{a}_{j,p}
\Biggr)\nonumber\\
&&\qquad\leq \frac{\sigma_{\varepsilon}^{2}}{n} \Biggl(k+ m + k + \sum_{l=m+k+1}^p
\mathbf{a}_{l,p}^{\tau} (\rho\mathbf{A}_{m})^{-{1}/{2}}\nonumber\\
&&\hspace*{57pt}\hspace*{91.2pt}{}\times
\Biggl(\sum_{j=k+1}^p\lambda_{x,j} u_ju_j^{\tau}\Biggr)
(\rho\mathbf{A}_{m})^{-{1}/{2}} \mathbf{a}_{l,p} \Biggr)
\\
&&\qquad\leq \frac{\sigma_{\varepsilon}^{2}}{n} \Biggl( m + 2k + \frac{1}{\mu
_{m+k+1}\cdot\rho}\sum_{j=k+1}^p \lambda_{x,j} \Biggr) \nonumber\\
&&\qquad\leq \frac{\sigma_{\varepsilon}^{2}}{n} ( m + 2k +CkC_0)\nonumber\\
&&\qquad
= \frac{\sigma_{\varepsilon}^{2}}{n} \bigl( m + \bigl[\rho^{-{1}/({2m+2q+1})}\bigr]\bigr)
(2+CC_0 ).\nonumber
\end{eqnarray}
This proves Relation (\ref{resthm12})
and completes the proof of Theorem \ref{thm1}.

\subsection[Proof of Theorem 2]{Proof of Theorem \protect\ref{thm2}}

With $\widehat{d}_{i} = \int_{I} \widehat{\alpha} (t) X_{i} (t) \,dt -
\frac{1}{p} \sum_{j=1}^{p}\widehat{\alpha} (t_{j}) X_{i} (t_{j}) $ we have
%
\begin{eqnarray}\label{L2proof00}
\| \widehat{\alpha} - \alpha\|_{\Gamma_n}^{2}
& \leq& \frac
{2}{n} \sum_{i=1}^{n} \Biggl[ \langle(X_{i}-\overline{X}) , \widehat
{\alpha} -
\alpha\rangle\nonumber\\
&&\hspace*{28.4pt}{} - \frac{1}{p} \sum_{j=1}^{p} \bigl(X_{i}(t_j)-\overline
{X}(t_{j})\bigr) \bigl(\widehat{\alpha} (t_{j}) - \alpha(t_{j})\bigr) \Biggr]^2
\nonumber\\[-8pt]\\[-8pt]
& &{} + \frac{2}{n} \sum_{i=1}^{n} \Biggl[ \frac{1}{p} \sum_{j=1}^{p}
(X_{i}-\overline{X}) (t_{j}) \bigl(\widehat{\alpha} (t_{j}) - \alpha(t_{j})\bigr)
\Biggr]^2 \nonumber\\
& \leq& \frac{4}{n} \sum_{i=1}^{n} ( \widehat{d}_{i}- \overline
{\widehat{d}})^{2} + \frac{4}{n} \sum_{i=1}^{n}
( d_{i}-\overline{d})^2 +
2 \| \widehat{\bolds{\alpha}} - \bolds{\alpha} \|
_{\Gamma_{n,p}}^{2}.
\nonumber
\end{eqnarray}
By assumptions (A.1)--(A.3), it follows from Theorem
\ref{thm1}, (\ref{H31}) and (\ref{H3a}) that the assertion of Theorem \ref{thm2} holds,
provided that
%
\begin{equation}\label{endthm2}
\frac{1}{n} \sum_{i=1}^{n} ( \widehat{d}_{i}- \overline{\widehat
{d}})^{2}
= O_P(p^{-2\kappa}).
\end{equation}
The proof of (\ref{endthm2}) consists of several steps. We will start
by giving a stochastic bound for
$\frac{1}{p}\widehat{\bolds{\alpha}}^\tau
\widehat{\bolds{\alpha}}$ and then study the stochastic behavior
of $\int_{0}^1 \widehat{\alpha}^{(m)} (t)^{2} \,dt$.
The use of a suitable Taylor expansion will then lead to the desired result.

By definition of $\widehat{\bolds{\alpha}}$ we have
\begin{eqnarray}\label{alphahbound}
\frac{1}{p}\widehat{\bolds{\alpha}}^\tau
\widehat{\bolds{\alpha}}
&\le&
\frac{3}{p} \bolds{\alpha}^\tau
\frac{1}{n p} \mathbf{X}^{\tau} \mathbf{X}
\biggl( \frac{1}{n p} \mathbf{X}^{\tau} \mathbf{X}
+ \rho\mathbf{A}_{m} \biggr)^{-2}\frac{1}{n p} \mathbf{X}^{\tau}
\mathbf{X}
\bolds{\alpha}
\nonumber\\[-2pt]
&&{} +3 \frac{1}{n^2 p}\mathbf{d}^\tau \mathbf{X}
\biggl( \frac{1}{n p} \mathbf{X}^{\tau} \mathbf{X}
+ \rho\mathbf{A}_{m} \biggr)^{-2} \mathbf{X}^{\tau}
\mathbf{d}\\[-2pt]
&&{}
+3 \frac{1}{n^2 p}\bolds{\varepsilon}^\tau \mathbf{X}
\biggl( \frac{1}{n p} \mathbf{X}^{\tau} \mathbf{X}
+ \rho\mathbf{A}_{m} \biggr)^{-2}\mathbf{X}^{\tau}
\bolds{\varepsilon}.
\nonumber
\end{eqnarray}
Since all eigenvalues of the matrix $\frac{1}{n p} \mathbf{X}^{\tau}
\mathbf{X}
( \frac{1}{n p} \mathbf{X}^{\tau} \mathbf{X}
+ \rho\mathbf{A}_{m} )^{-2}\frac{1}{n p} \mathbf{X}^{\tau}
\mathbf{X}$ are less than or equal to 1, the
first term on the right-hand side of (\ref{alphahbound}) is less than
or equal to $\frac{3}{p}\bolds{\alpha}^\tau
\bolds{\alpha}=O(1)$.
It is easily seen that the smallest eigenvalue of the matrix $\frac
{1}{n p}\mathbf{X}( \frac{1}{n p} \mathbf{X}^{\tau} \mathbf{X}+
\rho\mathbf{A}_{m} )^{-2} \mathbf{X}^{\tau}$ is proportional to
$1/\rho$, and thus the second term can be bounded by a term of order
$p^{-2\kappa}/\rho$. By (\ref{egams}) the expected value of the third
term is bounded by
\[
\frac{\sigma_{\varepsilon}^{2}}{n} \operatorname{Tr} \biggl[\frac{1}{n p}\mathbf{X}
\biggl( \frac{1}{n p} \mathbf{X}^{\tau} \mathbf{X}
+ \rho\mathbf{A}_{m} \biggr)^{-2} \mathbf{X}^{\tau}\biggr]\le\frac{\sigma
_{\varepsilon}^{2}}{n} \operatorname{Tr} [
( \rho\mathbf{A}_{m})^{-1}]=O\bigl(1/(n\rho)\bigr).
\]
We therefore arrive at
%
\begin{equation}\label{alphahboundnew}
\frac{1}{p}\widehat{\bolds{\alpha}}^\tau
\widehat{\bolds{\alpha}} =O_P\biggl(1+\frac{p^{-2\kappa}}{\rho
}+\frac{1}{n\rho}
\biggr).
\end{equation}
As a next step we will study the asymptotic behavior of $\int
_{0}^1 \widehat{\alpha}^{(m)} (t)^{2} \,dt$.
Since $\widehat{\bolds{\alpha}}$ is solution of the minimization
problem (\ref{nonnoisypbmino}), we can write
\begin{eqnarray*}
&&\frac{1}{n} \biggl\| \mathbf{Y} - \frac{1}{p} \mathbf{X}
\widehat{\bolds{\alpha}} \biggr\|^{2} + \frac{\rho}{p}\widehat
{\bolds{\alpha}}^\tau\mathbf{P}_m\widehat{\bolds{\alpha
}}+ \rho\int_{0}^1 \widehat{\alpha}^{(m)} (t)^{2} \,dt \\[-2pt]
&&\qquad\leq\frac{1}{n} \biggl\| \mathbf{Y} - \frac{1}{p} \mathbf{X}
\bolds{\alpha} \biggr\|^{2} + \frac{\rho}{p}\bolds{\alpha
}^\tau\mathbf{P}_m\bolds{\alpha} + \rho\int_{0}^1 \alpha
^{(m)} (t)^{2} \,dt ,
\end{eqnarray*}
and therefore
%
\begin{eqnarray}\label{intalphammin}
 \rho\int_{0}^1 \widehat{\alpha}^{(m)} (t)^{2} \,dt
 & \leq&
\| \widehat{\bolds{\alpha}} -
\bolds{\alpha} \|_{\Gamma_{n,p}}^{2} + \frac{2}{n} \biggl\langle
\mathbf{Y} - \frac{1}{p} \mathbf{X} \bolds{\alpha} , \frac{1}{p}
\mathbf{X} \widehat{\bolds{\alpha}} - \frac{1}{p} \mathbf{X}
\bolds{\alpha} \biggr\rangle\nonumber\\[-9pt]\\[-9pt]
&&{} + \rho\int_{0}^1 \alpha^{(m)} (t)^{2} \,dt - \frac{\rho}{p}\widehat
{\bolds{\alpha}}^\tau\mathbf{P}_m\widehat{\bolds{\alpha}}
+ \frac{\rho}{p}\bolds{\alpha}^\tau\mathbf{P}_m\bolds
{\alpha}. \nonumber
\end{eqnarray}
We have to focus on the term
\[
\frac{2}{n} \biggl\langle\mathbf{Y} - \frac{1}{p} \mathbf{X} \bolds
{\alpha} ,
\frac{1}{p} \mathbf{X} \widehat{\bolds{\alpha}} - \frac{1}{p}
\mathbf{X} \bolds{\alpha} \biggr\rangle
=\frac{2}{n} \biggl\langle\mathbf{d}+ \bolds{\varepsilon},
\frac{1}{p} \mathbf{X} \widehat{\bolds{\alpha}} - \frac{1}{p}
\mathbf{X} \bolds{\alpha} \biggr\rangle.
\]
The Cauchy--Schwarz inequality together with the definition of $\| \cdot
\|_{\Gamma_{n,p}}^{2}$ yield
%
\begin{equation}\label{intalnew1}
\frac{1}{n} \mathbf{d}^\tau\biggl(\frac{1}{p} \mathbf{X} \widehat
{\bolds{\alpha}} - \frac{1}{p} \mathbf{X} \bolds{\alpha
}\biggr)=O_P(p^{-\kappa} \| \widehat{\bolds{\alpha}} - \bolds
{\alpha} \|_{\Gamma_{n,p}}).
\end{equation}
Note that
\begin{eqnarray*}
&&\frac{2}{n} \biggl\langle \bolds{\varepsilon},
\frac{1}{p} \mathbf{X} \widehat{\bolds{\alpha}} - \frac{1}{p}
\mathbf{X} \bolds{\alpha} \biggr\rangle\\
&&\qquad=\frac{2}{n}\bolds{\varepsilon}^\tau\biggl(\frac{1}{p}\mathbf{X}
\mathbb{E}_\varepsilon( \widehat{\bolds{\alpha}}) - \frac{1}{p}
\mathbf{X} \bolds{\alpha}\biggr)+
\frac{2}{n}\bolds{\varepsilon}^\tau\biggl(\frac{1}{p}\mathbf{X}
\widehat{\bolds{\alpha}} - \frac{1}{p}\mathbf{X}\mathbb
{E}_\varepsilon( \widehat{\bolds{\alpha}})\biggr).
\end{eqnarray*}
Obviously, $\frac{1}{n}\bolds{\varepsilon}^\tau(\frac{1}{p}\mathbf{X}
\mathbb{E}_\varepsilon( \widehat{\bolds{\alpha}})-\frac{1}{p} \mathbf
{X} \bolds{\alpha})$ is a zero mean random variable with variance
bounded by
$\frac{\sigma^2_\varepsilon}{n}\|
\mathbb{E}_\varepsilon( \widehat{\bolds{\alpha}})-\bolds{\alpha
}\|_{\Gamma_{n,p}}^{2}$. By definition of
$\widehat{\bolds{\alpha}}$, (\ref{H31}), (\ref{thm1step1}) and
(\ref{thm1step4}) we have
\begin{eqnarray*}
\mathbb{E}_\varepsilon\biggl(\frac{1}{n}\bolds{\varepsilon}^\tau\biggl(\frac
{1}{p}\mathbf{X}
\widehat{\bolds{\alpha}} - \frac{1}{p}\mathbf{X}\mathbb
{E}_\varepsilon( \widehat{\bolds{\alpha}})\biggr)\biggr)
& \le&
\frac{\sigma_{\varepsilon}^{2}}{n} \operatorname{Tr} \biggl[ \biggl( \frac
{1}{n p} \mathbf{X}^{\tau} \mathbf{X}
+ \rho\mathbf{A}_{m} \biggr)^{-1} \frac{1}{n p} \mathbf{X}^{\tau}
\mathbf{X} \biggr] \\
& = & O_P\biggl(\frac{1}{n\rho^{{1}/({2m+2q+1})}}\biggr).
\end{eqnarray*}
We can conclude that
%
\begin{equation}\label{intalnew2}\qquad \ \
\frac{2}{n} \biggl\langle \bolds{\varepsilon},
\frac{1}{p} \mathbf{X} \widehat{\bolds{\alpha}} - \frac{1}{p}
\mathbf{X} \bolds{\alpha} \biggr\rangle=
O_P\biggl(\frac{1}{\sqrt{n}}\|
\mathbb{E}_\varepsilon( \widehat{\bolds{\alpha}}-\bolds{\alpha
})\|_{\Gamma_{n,p}}+\frac{1}{n\rho^{{1}/({2m+2q+1})}}\biggr).
\end{equation}
When combining (\ref{alphahboundnew}), (\ref{intalphammin}), (\ref
{intalnew1}) and (\ref{intalnew2}) with the results of Theorem \ref{thm1} we
thus obtain
%
\begin{equation}\label{intalbound}
\int_{0}^1 \widehat{\alpha}^{(m)} (t)^{2} \,dt=O_P\biggl(1+\frac
{p^{-2\kappa}}{\rho}+\frac{1}{n\rho^{({2m+2q+2})/({2m+2q+1})}}
\biggr).
\end{equation}
Let us now expand $\widehat{\alpha}$ into a Taylor series:
$\widehat{\alpha} (t) = P(t) + R(t)$ for all
$t \in[0,1]$ with
\begin{eqnarray*}
P(t)
&=& \sum_{l=0}^{m-1} \frac{t^{l}}{l!} \widehat{\alpha}^{(l)} (0),\qquad
R(t)=\int_0^t r(s)\,ds
\end{eqnarray*}
and
\begin{eqnarray*}
r(t) &=& \int_{0}^{t} \frac{(t-u)^{m-1}}{(m-1)!} \widehat{\alpha}^{(m)}
(u) \,du.
\end{eqnarray*}
It follows from (\ref{alphahboundnew}) as well as (\ref
{intalbound}) that
$|\widehat{\alpha}^{(l)} (0)|=O_P(1+(\frac{p^{-2\kappa}}{\rho
})^{1/2}+(\frac{1}{n\rho^{({2m+2q+2})/({2m+2q+1})}})^{1/2})$
for $l=0,\ldots,m-1$, and some straightforward calculations
yield
\begin{eqnarray*}
\biggl|  \| \widehat{\alpha} \|^2 -\frac{1}{p}\widehat{\bolds
{\alpha}}^\tau
\widehat{\bolds{\alpha}}\biggr| &=& \biggl|\int_0^1
\bigl(P(t)+R(t)\bigr)^2\,dt-\frac{1}{p}\sum_{j=1}^p \bigl( P(t_j)+R(t_j)\bigr)^2\biggr|
\\
&\le&
\Biggl(\sum_{j=1}^p \biggl[\int
_{t_j-1/(2p)}^{t_j+1/(2p)}\bigl(P(t)+R(t)+P(t_j)+R(t_j)\bigr)^2\,dt\biggr]^2
\Biggr)^{1/2} \\
&&{}\times\Biggl(\sum_{j=1}^p\frac{1}{p}\biggl[\int
_{t_j-1/(2p)}^{t_j+1/(2p)}|P'(s)|+ |r(s)|\,ds\biggr]^2 \Biggr)^{1/2},
\end{eqnarray*}
which leads to
%
\begin{eqnarray}\label{l2proof0}
&&\biggl|  \| \widehat{\alpha} \|^2 -\frac{1}{p}\widehat{\bolds
{\alpha}}^\tau
\widehat{\bolds{\alpha}}\biggr|\nonumber\\[-8pt]\\[-8pt]
&&\qquad=O_P\biggl(p^{-1}\cdot\biggl(1+\frac
{p^{-2\kappa}}{\rho}+
\bigl[n\rho^{({2m+2q+2})/({2m+2q+1})}\bigr]^{-1}\biggr)\biggr).\nonumber
\end{eqnarray}
Using again (\ref{alphahboundnew}) and our assumptions on
$\rho,p,n$, this implies
%
\begin{equation}\label{gammaalpha0}
\| \widehat{\alpha} \|^{2} =O_P(1).
\end{equation}
At the same time, (\ref{alphahboundnew}) and (\ref{intalbound})
together with assumptions (A.1) and (A.2) imply that
with $\widetilde{X}_i=X_i-\overline{X}$
\begin{eqnarray*}
\frac{1}{n}\sum_{i=1}^n(\widehat{d}_{i}-\overline{\widehat
{d}})^2
&=&
\frac{1}{n}\sum_{i=1}^n\Biggl( \sum_{j=1}^p\int
_{t_j-1/(2p)}^{t_j+1/(2p)}\bigl(\widehat{\alpha} (t)-
\widehat{\alpha} (t_j)\bigr)\widetilde{X}_i(t)\nonumber\\
&&\hspace*{61pt}\hspace*{32.7pt}{} +\widehat{\alpha}
(t_j)\bigl(\widetilde{X}_i(t)- \widetilde{X}_i(t_j)\bigr)\,dt\Biggr)^2 \nonumber\\
&\le&
2x_{\max}^2 \Biggl(
\sum_{j=1}^p \frac{1}{p} \biggl[\int_{t_j-1/(2p)}^{t_j+1/(2p)} |P'(t)|+
|r(t)|\,dt\biggr]^2\Biggr)
\nonumber\\
&&{} +2 \Biggl(\frac{1}{p} \sum_{j=1}^{p}
\widehat{\alpha} (t_{j})^2\Biggr)
\frac{1}{n}\sum_{i=1}^n
\sum_{j=1}^p \int_{t_j-1/(2p)}^{t_j+1/(2p)}\bigl(\widetilde
{X}_i(t)-\widetilde{X}(t_j)\bigr)^2\,dt\nonumber
\end{eqnarray*}
and thus
%
%
\begin{eqnarray}\label{L2proof01}\qquad
&&\frac{1}{n}\sum_{i=1}^n(\widehat{d}_{i}-\overline{\widehat{d}})^2 =
O_P\biggl(p^{-2}\biggl(1+\frac{p^{-2\kappa}}{\rho}+\frac{1}{n\rho^{(
{2m+2q+2})/({2m+2q+1})}}\biggr)\nonumber\\[-8pt]\\[-8pt]
&&\hspace*{172pt}{}
+p^{-2\kappa}\biggl(1+\frac{p^{-2\kappa}}{\rho}+\frac{1}{n\rho}\biggr)\biggr).\nonumber
\end{eqnarray}
By our assumptions on $\rho,p,n$, relation (\ref{endthm2}) is an
immediate consequence. This completes the proof
of Theorem \ref{thm2}.

\subsection[Proof of Theorem 3]{Proof of Theorem \protect\ref{thm3}}

In terms of eigenvalues and eigenfunctions of $\Gamma$ we obviously obtain
\[
\langle\Gamma u , u\rangle=\sum_r \lambda_r \langle\zeta_r , u\rangle^2.
\]
Let $\tau_{ri}=\langle X_i-\mathbb{E}(X),\zeta_r\rangle$ for $r=1,2,\ldots$
and $i=1,\ldots,n$. Some well-known results of stochastic process
theory now can be summarized as follows:
\begin{longlist}
\item $\mathbb{E}(\tau_{ri})=0$, $\mathbb{E}(\tau_{ri}^2)=\lambda
_r$, and
$\mathbb{E}(\tau_{ri}\tau_{si})=0$ for all $r,s$, $s\neq r$ and
$i=1,\ldots,n$.
\item For any $k=1,2,\ldots,$ the eigenfunctions $\zeta_1,\ldots
,\zeta_k$ corresponding to
$\lambda_1\ge\cdots\ge\lambda_k$ provide a best basis for approximating
$X_i$ by a $k$-dimensional linear space:
\end{longlist}
%
%
\begin{eqnarray}\label{gammaalphasub}
\sum_{r=q+1}^\infty\lambda_r &=&\mathbb{E}\Biggl(\Biggl\| X-\mathbb
{E}(X)-\sum_{s=1}^q \langle
X-\mathbb{E}(X),\zeta_s\rangle\zeta_s\Biggr\|^2\Biggr)\nonumber\\[-8pt]\\[-8pt]
&\le&\mathbb{E}\biggl( \inf_{f\in\mathcal{L}_k} \| X-\mathbb{E}(X) -f\|
^2\biggr),\nonumber
\end{eqnarray}
for any other $k$-dimensional linear subspace $\mathcal{L}_k$ of $L^2([0,1])$.

By (A.3) we can conclude that
%
\begin{equation}\label{gammalambdarate}
\sum_{r=k+1}^\infty\lambda_r =O(k^{-2q}) \qquad\mbox{as } k\rightarrow
\infty.
\end{equation}
At first we have
\[
\| \widehat{\alpha} - \alpha\|_{\Gamma_n}^{2} \le\frac
{2}{n}\sum_{i=1}^n\langle\widehat{\alpha} - \alpha, X_i-\mathbb
{E}(X)\rangle^2+\frac{2}{n}\sum_{i=1}^n\langle\widehat{\alpha} - \alpha
, \mathbb{E}(X)-\overline{X}\rangle^2,
\]
and by (\ref{gammaalpha0}) and with assumption (A.4) the last term is
of order $O_P(n^{-1})$. The relevant semi-norms can now be rewritten in
the form
%
\begin{equation}\label{gammaalphadef}
\| \widehat{\alpha} - \alpha\|_{\Gamma}^{2}=
\sum_{r=1}^\infty\lambda_r \langle\zeta_r , \widehat{\alpha} - \alpha
\rangle^2=:
\sum_{r=1}^\infty\lambda_r \widetilde\alpha_r^2
\end{equation}
and
%
\begin{eqnarray}\label{gammaalphader}\quad
\| \widehat{\alpha} - \alpha\|_{\Gamma_n}^{2}&=&
\| \widehat{\alpha} - \alpha\|_{\Gamma}^{2}+
\sum_{r=1}^\infty\sum_{s=1}^\infty \widetilde\alpha_r\widetilde\alpha
_s \Biggl(
\frac{1}{n}\sum_{i=1}^n\tau_{ri}\tau_{si}-\lambda_rI(r=s)\Biggr)\nonumber\\[-8pt]\\[-8pt]
&&{} + O_P(n^{-1}),\nonumber
\end{eqnarray}
where $I(r=s)=1$ if $r=s$, and $I(r=s)=0$ if $r\ne s$.
Define
\[
\widetilde{\tau}_{rr}=
\frac{1}{\lambda_r\sqrt{n}} \sum_{i=1}^n (\tau_{ri}^2-\lambda_r)
\quad\mbox{and}\quad
\widetilde{\tau}_{rs}=
\frac{1}{\sqrt{\lambda_r\lambda_sn}} \sum_{i=1}^n \tau_{ri}\tau_{si},\  r\neq s
\]
(with $\widetilde{\tau}_{rs}:=0$ if $\min\{\lambda_r,\lambda_s\}=0$).
The properties of $\tau_{ri}$ given in (i) imply that $\mathbb
{E}(\widetilde{\tau}_{rs})=0$ for all $r,s$, and
we can infer from assumption (A.4) that for some $C_{10}<\infty$
%
\begin{equation}\label{gammaalphabound}
\mathbb{E}(\widetilde{\tau}_{rs}^2)\le C_{10},
\end{equation}
holds for all $r,s=1,2,\ldots$ and all sufficiently large $n$.
Using the Cauchy--Schwarz inequality we therefore obtain for all
$k=0,1,\ldots$
\begin{eqnarray}\label{gammaalphacauchy}
&&\Biggl|\sum_{r=1}^\infty\sum_{s=1}^\infty \widetilde\alpha_r\widetilde
\alpha_s \Biggl(
\frac{1}{n}\sum_{i=1}^n\tau_{ri}\tau_{si}-\lambda_rI(r=s)\Biggr)\Biggr|\nonumber\\
&&\qquad= \Biggl|\frac{1}{\sqrt{n}}
\sum_{r=1}^\infty\sum_{s=1}^\infty \widetilde\alpha_r\widetilde\alpha_s
(\lambda_r\lambda_s)^{{1}/{2}}\widetilde{\tau}_{rs}\Biggr|
\nonumber\\[-8pt]\\[-8pt]
&&\qquad  \le
\frac{2}{\sqrt{n}}\Biggl(\sum_{r=1}^k\sum_{s= r}^\infty\lambda_r
\widetilde\alpha_r^2\widetilde\alpha_s^2\Biggr)^{{1}/{2}}
\Biggl(\sum_{r=1}^k\sum_{s= r}^\infty\lambda_s \widetilde{\tau}_{rs}^2
\Biggr)^{{1}/{2}} \nonumber\\
&&\qquad\quad{} +
\frac{2}{\sqrt{n}}\Biggl(\sum_{r=k+1}^\infty\sum_{s= r}^\infty
\widetilde\alpha_r^2\widetilde\alpha_s^2\Biggr)^{{1}/{2}}
\Biggl(\sum_{r=k+1}^\infty\sum_{s= r}^\infty\lambda_r\lambda_s
\widetilde{\tau}_{rs}^2 \Biggr)^{{1}/{2}}. \nonumber
\end{eqnarray}
Relation (\ref{gammaalpha0}) leads to $\| \widehat{\alpha} -\alpha\|^2
\ge\sum_{r=1}^\infty\widetilde\alpha_r^2=O_P(1)$,
which together with (\ref{gammaalphadef}) implies that for arbitrary $k$
\[
\Biggl(\sum_{r=1}^k\sum_{s= r}^\infty\lambda_r \widetilde\alpha_r^2\widetilde
\alpha_s^2\Biggr)^{{1}/{2}}
\le\Biggl(\Biggl(\sum_{r=1}^\infty\lambda_r \widetilde\alpha_r^2\Biggr)\Biggl(\sum_{s=
1}^\infty\widetilde\alpha_s^2\Biggr)\Biggr)^{{1}/{2}}
=O_P(\| \widehat{\alpha} - \alpha\|_{\Gamma}).
\]
Choose $k$ proportional to $n^{1/2}$. Relation (\ref{gammalambdarate})
then yields
$\sum_{r=k+1}^\infty\sum_{s= r}^\infty\lambda_r\lambda_s \le(\sum
_{r=k+1}^\infty\lambda_r)^2
=O(n^{-2q})$ and $\sum_{r=1}^k\sum_{s= r}^\infty\lambda_s =O(\max\{
\log n,n^{(1-2q)/2}\})$.
Since by (\ref{gammaalphabound}) the moments of $\widetilde{\tau}_{rs}$
are uniformly bounded for all $r,s$,
it follows that
\begin{eqnarray*}
\Biggl(\sum_{r=1}^k\sum_{s= r}^\infty\lambda_s
\widetilde{\tau}_{rs}^2\Biggr)^{{1}/{2}}
&=&O_P\bigl(\max\bigl\{\log n,n^{(1-2q)/4}\bigr\}\bigr)
,\\
\Biggl(\sum_{r=k+1}^\infty\sum_{s= r}^\infty\lambda_r\lambda_s \widetilde
{\tau}_{rs}^2 \Biggr)^{{1}/{2}}&=&O_P(n^{-q}).
\end{eqnarray*}
When combining these results we can conclude that
\begin{eqnarray*}
&& %
\Biggl|\sum_{r=1}^\infty\sum_{s=1}^\infty \widetilde\alpha_r\widetilde
\alpha_s \Biggl(
\frac{1}{n}\sum_{i=1}^n\tau_{ri}\tau_{si}-\lambda_rI(r=s)\Biggr)\Biggr| \\
&&\qquad =O_P\bigl(
\max\bigl\{n^{-1/2}\log n\cdot\| \widehat{\alpha} - \alpha\|_{\Gamma},
n^{-(2q+1)/4}\cdot\| \widehat{\alpha} - \alpha\|_{\Gamma},
n^{-(2q+1)/2}\bigr\}\bigr).
\end{eqnarray*}
Together with (\ref{gammaalphader})
assertion (\ref{resthm3}) now follows from the rates of convergence of
$\| \widehat{\alpha} - \alpha\|_{\Gamma_n}^{2}$ derived in Theorem \ref{thm2}.

It remains to prove (\ref{prederror}). Note that by our assumptions on
$\varepsilon_i$ and assumption~(A.4) we have
$|\mathbb{E}(Y)-\overline{Y}|^2\le2 \overline{\varepsilon}^2+2 \langle
\alpha,\mathbb{E}(X)-\overline{X}\rangle^2=O_P(n^{-1})$.
Together with (\ref{gammaalpha0}) and assumption (A.4) this implies
%
\begin{eqnarray*}
&&\bigl|\mathbb{E} \bigl((\widehat{\alpha_0}+\langle\widehat{\alpha
},X_{n+1}\rangle-
\alpha_0-\langle\alpha,X_{n+1}\rangle
)^2| \widehat{\alpha_0},\widehat{\alpha}
\bigr)-\| \widehat{\alpha} - \alpha\|_{\Gamma}^{2}\bigr|\\
&&\qquad\le2|\mathbb{E}(Y)-\overline{Y}|^2+2 \langle\widehat{\alpha},\mathbb
{E}(X)-\overline{X}\rangle^2
= O_P(n^{-1}),
\end{eqnarray*}
which completes the proof of the theorem.

\subsection[Proof of Proposition 1]{Proof of Proposition \protect\ref{prop1}}

In dependence of $q$ we first construct special probability
distributions of $X_i$. For $2q=1$,
$\tau\in[0,1]$ and $r:=0$ set $\widetilde X_{\tau;0}(t):=1$ for $t\in
[0,\tau]$ and
$\widetilde X_{\tau,0}(t):=0$
for $t\in(\tau,1]$. For $2q\ge3$, $\tau\in[0,1]$, and $r:=q-0.5$ let
$\widetilde
X_{\tau;r}(t):=\frac{1}{r!}t^{r}$ for $t\in[0,\tau]$ and $\widetilde
X_{\tau;r}(t):=\sum_{j=0}^{r-1}
\frac{1}{(r-j)!} \tau^{r-j}(t-\tau)^j$ for $t\in(\tau,1]$.

For $k=1,2,\ldots$ let $\mathcal{L}_{(r+1)k}$ denote the $(r+1)\cdot k$
dimensional linear space of all
functions $g_\beta$ of the form $g_\beta(t):=\sum_{j=0}^{k-1} (\sum
_{l=0}^r \beta_{l,j} t^l)\cdot I(t\in
[\frac{j}{k},\frac{j+1}{k}])$. It is then easily verified that $\sup
_{t\in[{j}/{k},({j+1})/{k}]}
\min_\beta|g_\beta(t)-\widetilde X_{\tau;r}(t)|=0$ if $\tau\notin
[\frac{j}{k},\frac{j+1}{k}]$, while
$\sup_{t\in[{j}/{k},({j+1})/{k}]} \min_\beta|g_\beta
(t)-\widetilde X_{\tau;r}(t)|\le k^{-r}$
if $\tau\in[\frac{j}{k},\frac{j+1}{k}]$. It follows that there exist
constants $B_r\le1$ such that the functions $B_r\widetilde X_{\tau;r}(t)$
satisfy $\inf_{g_\beta\in\mathcal{L}_{(r+1)k}} \int_0^1 (B_r\widetilde
X_{\tau;r}(t)-g_\beta(t))^2\,dt \le C
(r+2)^{-(2r+1)}k^{-(2r+1)}=C(r+2)^{-2q}k^{-2q}$ for all $k=1,2,\ldots.$

Now let $\tau_1,\ldots,\tau_n$ denote i.i.d. real random variables which
are uniformly distributed on
$[0,1]$ and let $X_{\tau_i;r}=B_r\widetilde X_{\tau_i;r}(t)-\mathbb
{E}(B_r\widetilde X_{\tau_i;r}(t))$.
Obviously, $\tau_i\rightarrow X_{\tau_i,r}^{(j)}(t)$ is a continuous
mapping from $[0,1]$ on
$L^2([0,1])$, and
the
probability distribution of $\tau_i$ induces a corresponding centered
probability distribution $P_{r}$ on
$L^2([0,1])$. Since the eigenfunctions of the corresponding covariance
operator provide a best basis for
approximating $X_i$ by a $k$-dimensional linear space, we obtain from
what is done above
\[
\sum_{j=k+1}^\infty\lambda_j \le\mathbb{E}\biggl( \inf_{g_\beta^*\in
\mathcal{L}^*_{(r+1)[{k}/({r+1})]}} \|X_{\tau_i;r} -g_\beta^*\|^2\biggr)
\le Ck^{-2q},
\]
for all sufficiently large $k$ and $\mathcal{L}_{(r+1)k}^*:=\{g_\beta
-\mathbb{E}(B_r\widetilde X_{\tau_i;r})|g_\beta
\in\mathcal{L}_{(r+1)k}\}$.

In order to verify that $P_{r}\in\mathcal{P}_{q,C}$, it remains to check
the behavior of
$\|\widehat{\alpha}-\alpha\|_{\Gamma}=\int_0^1 \langle X_{\tau
;r},\widehat{\alpha}-\alpha
\rangle^2 \,d\tau$.
First
note that although assumption (A.2) does not hold for $2q=1$, even in
this case, with $\kappa=1/2$,
relation (\ref{H3a}) holds and arguments in the proof of Theorems \ref{thm1}
and \ref{thm2} imply that for sufficiently large $p$,
$\frac{1}{n} \sum_{i=1}^n \langle X_{\tau_i;r},\widehat{\alpha}-\alpha
\rangle^2= O_P(n^{{-(2m+2q+1)}/({2m+2q+2})})$. For some $1>\delta
>\frac{2m+2q+1}{2m+2q+2}$ define a partition of
$[0,1]$ into $n^{\delta}$ disjoint intervals $I_1,\ldots,I_{n^\delta}$
of equal length $n^{-\delta}$.
For $j=1,\ldots, n^{\delta}$, let
$s_j$ denote the midpoint of the interval $I_j$, and use $n_j$ denote
the (random) number of $\tau_1,\ldots,\tau_n$
falling into $I_j$. By using the Cauchy--Schwarz inequality as well as a
definition of $X_{\tau;r}$ it is easily verified
that there exists a constant $L_r<\infty$ such that $|\langle X_{\tau
;r},\widehat{\alpha}-\alpha\rangle-
\langle X_{\tau^*;r},\widehat{\alpha}-\alpha\rangle|\le L_r |\tau-\tau
^*|^{1/2}\| \widehat{\alpha}-\alpha\|$
for $\tau,\tau^*\in[0,1]$
($|\tau-\tau^*|^{1/2}$ may be replaced by $|\tau-\tau^*|$ if $2q>1$). Then
\begin{eqnarray*}
&&|\langle X_{\tau;r},\widehat{\alpha}-\alpha\rangle^2-
\langle X_{\tau^*;r},\widehat{\alpha}-\alpha\rangle^2|\\
&&\qquad \le2L_r |\tau-\tau^*|^{1/2}\| \widehat{\alpha}-\alpha\|
\min\{|\langle X_{\tau;r},\widehat{\alpha}-\alpha\rangle|,|\langle
X_{\tau^*;r},\widehat{\alpha}-\alpha\rangle|\}\\
&&\qquad\quad{}
+L_r^2 |\tau-\tau^*|\| \widehat{\alpha}-\alpha\|^2.
\end{eqnarray*}
By (\ref{gammaalpha0}) another application of the Cauchy--Schwarz inequality
leads to\break $\frac{1}{n} \sum_{i=1}^n \langle X_{\tau_i;r},\widehat{\alpha
}-\alpha
\rangle^2=\frac{1}{n} \sum_{j=1}^{n^{\delta}}
n_j\langle X_{s_j;r},\widehat{\alpha}-\alpha
\rangle^2+
o_P(n^{{-(2m+2q+1)}/({2m+2q+2})})$. Since $\sup_{j=1,\ldots,n^\delta}
\frac{|n_j-\mathbb{E}( n_j)|}{n_j}=O_P(1)$
with
$\mathbb{E}( n_j)=n\cdot n^{-\delta}$, we can conclude that
$\frac{1}{n} \sum_{j=1}^{n^{\delta}}
\mathbb{E}( n_j) \langle X_{s_j;r},\widehat{\alpha}-\alpha
\rangle^2=O_P(n^{{-(2m+2q+1)}/({2m+2q+2})})$. Finally,
%
\begin{eqnarray*}
&&\Biggl|\int_0^1 \langle X_{\tau;r},\widehat{\alpha}-\alpha
\rangle^2 \,d\tau-\frac{1}{n} \sum_{j=1}^{n^{\delta}}
\mathbb{E}( n_j) \langle X_{s_j;r},\widehat{\alpha}-\alpha
\rangle^2\Biggr|\\
&&\qquad\le\frac{1}{n^\delta} \sum_{j=1}^{n^{\delta}}\sup_{\tau\in I_j}
|\langle X_{\tau;r},\widehat{\alpha}-\alpha\rangle^2-
\langle X_{s_j;r},\widehat{\alpha}-\alpha\rangle^2|\\
&&\qquad=o_P\bigl(n^{
{-(2m+2q+1)}/({2m+2q+2})}\bigr),
\end{eqnarray*}
%
and the desired result $\|\widehat{\alpha}-\alpha\|_{\Gamma}=
O_P(n^{{-(2m+2q+1)}/({2m+2q+2})})$
is an immediate consequence. Therefore, $P_{r}\in\mathcal{P}_{q,C}$.

We now have to consider the functionals $\langle X_{\tau_i;r},\alpha
\rangle$ more closely. Let $\mathcal{C}^*(m+r+1,D)$
denote the space of all $m+r+1$-times continuously differentiable
functions $\widetilde\alpha$ satisfying
$\int_0^1 \widetilde\alpha(t)\,dt=0$ as well as
$\int_0^1 \widetilde\alpha^{(j)}(t)^2\,dt\le D$ for all $j=0,1,\ldots
,m+r+1$ as well as
$\widetilde\alpha^{(j)}(0)=\widetilde\alpha^{(j)}(1)=0$ for all
$j=0,\ldots,r+1$, and set $\mathcal{C}^*(m,r,D)=\{\alpha| \alpha=\widetilde\alpha^{(r+1)}, \widetilde
\alpha\in\mathcal{C}^*(m+r+1,D)\}$. Then, for any $\alpha\in
\mathcal{C}^*(m,0,D)$ there is a $\widetilde\alpha\in\mathcal{C}^*(m+1,D)$
such that
\begin{eqnarray*}
\langle X_{\tau_i;0},\alpha\rangle
&=&B_0\int_0^{\tau_i} \alpha(t)\,dt -\langle\mathbb{E}(B_0\widetilde
X_{\tau_i;0}),\alpha\rangle\\
&=&B_0 \widetilde\alpha(\tau_i)-B_0\int_0^{1} \widetilde\alpha
(t)\,dt=B_0 \widetilde\alpha(\tau_i)
\end{eqnarray*}
while for any $\alpha\in\mathcal{C}^*(m,r,D)$, $r\ge1$ and $\widetilde
\alpha\in\mathcal{C}^*(m+r+1,D)$,
$\alpha=\widetilde\alpha^{(r+1)}$, partial integration leads to
\begin{eqnarray*}
&&\langle X_{\tau_i;r},\alpha\rangle=(-1)^{r-1}\bigl\langle X_{\tau
_i;r}^{(r-1)},\widetilde\alpha^{(2)}\bigr\rangle
\\
&&\qquad =\bigl(X_{\tau_i;r}^{(r-1)}(\tau_i)\widetilde\alpha^{(2)}(\tau_i)-X_{\tau
_i;r}^{(r-1)}(0)\widetilde\alpha^{(2)}(0)\bigr)\\
&&\qquad\quad{}  +B_r(-1)^{r}\int_0^{\tau_i} \widetilde\alpha^{(1)}(t)\,dt
-B_r(-1)^{r}\mathbb{E}\biggl(\int_0^{\tau_i} \widetilde\alpha
^{(1)}(t)\,dt\biggr)\\
&&\qquad\quad{}  +
\bigl(X_{\tau_i;r}^{(r-1)}(1)\widetilde\alpha^{(2)}(1)-X_{\tau
_i;r}^{(r-1)}(\tau_i)\widetilde\alpha^{(2)}(\tau_i)\bigr)\\
&&\qquad =B_r(-1)^{r}\widetilde\alpha(\tau_i)-\mathbb{E}(B_r(-1)^{r}\widetilde
\alpha(\tau_i))=B_r(-1)^{r}\widetilde\alpha(\tau_i).
\end{eqnarray*}

Obviously, $\widetilde\alpha^*=B_r(-1)^{r}\widetilde\alpha\in\mathcal{C}^*(m+r+1,B_rD)$.
By construction,
with $f_{a}(\tau_i) :=\langle X_{\tau_i,r},a\rangle$ we generally obtain
\[
\| \alpha-\hat a(\alpha,P_{\beta})\|^2_\Gamma=\int_0^1
\bigl(f_{\alpha}(\tau)-f_{\hat a(\alpha,P_r)}(\tau)\bigr)^2\,d\tau.
\]
By definition, $f_{\alpha}(\tau)=\widetilde\alpha^*(\tau)=\mathbb
{E}(Y_i| \tau_i=\tau)$ is the regression
function in the regression model $Y_i=\widetilde\alpha^*(\tau_i)+
\varepsilon_i$, and we will use the notation
$S_n(\widetilde\alpha^*)$ to denote an estimator of $\widetilde\alpha
^*$ from the data
$(Y_i,\tau_i),\ldots,(Y_n,\tau_n)$. Note that knowledge of $(Y_i,\tau
_i)$ is equivalent to knowledge of
$(Y_i,X_{\tau_i;r})$, and an estimator $f_{\hat a(\alpha,P_r)}$ of
$\widetilde\alpha^*$ can thus be seen as a
particular estimator $S_n(\widetilde\alpha^*)$ based on $(Y_i,\tau
_i),\ldots,(Y_n,\tau_n)$. We can conclude
that as $n\rightarrow\infty$,
\begin{eqnarray*}
&&\sup\limits_{ P\in\mathcal{P}_{q,C}} \sup\limits_{\alpha\in\mathcal{C}_{m,D}} \inf\limits_{\hat
a(\alpha,P)}\mathbb{P}\bigl( \| \alpha- \hat a(\alpha,P)\|^2_\Gamma\\
&&\qquad\ge
c_n\cdot
n^{{-(2m+2q+1)}/({2m+2q+2})}
\bigr)\\
&&\qquad\ge\sup\limits_{\widetilde\alpha^*\in C^*(m+r+1,B_rD)} \inf\limits
_{S_n(\widetilde\alpha^*)}\mathbb{P}\biggl(\int_0^1
\bigl(\widetilde\alpha^*(\tau)-S_n(\widetilde\alpha^*)(\tau)
\bigr)^2\,d\tau\\
&&\hspace*{150.5pt}
\ge c_n\cdot n^{{-(2m+2q+1)}/({2m+2q+2})}\biggr)\rightarrow1.
\end{eqnarray*}
Convergence of the last probability to 1 follows from well-known
results on optimal
rates of convergence in nonparametric regression (cf. Stone \cite{St82}).

\subsection[Proof of Proposition 2]{Proof of Proposition \protect\ref{prop2}}

We first consider (\ref{gcvrho}).
The set $\{ \mathbf{H}_\rho\}_{\rho>0}$ constitutes an ordered linear
smoother according to the definition in Kneip~\cite{Kn94}. Theorem~1 of
Kneip \cite{Kn94}
then implies that
$|\mathit{MSE}_m(\hat\rho^*)-\mathit{MSE}_m(\rho_{\mathrm{opt}})|=O_P(n^{-1/2}\times
\mathit{MSE}_m(\rho_{\mathrm{opt}})^{1/2})$,
where $\hat\rho^*$ is determined by minimizing Mallow's $C_L$,
$C_L(\rho):=
\frac{1}{n}\|\mathbf{Y}-\mathbf{H}_\rho\mathbf{Y}\|^2+\frac{2\sigma
^2_\varepsilon}{n}\operatorname{Tr}(\mathbf{H}_\rho)$. Note that although we consider
centered values $Y_i-\overline{Y}$ instead of $Y_i$ all arguments in
Kneip \cite{Kn94} apply, since
$(\overline{Y},\ldots,\overline{Y})^\tau\mathbf{X}=0$.
The arguments used in the proof of Theorem 1 of Kneip (\cite{Kn94},
relations (A.17)--(A.22)) imply that for all $\rho$ the difference
$C_L(\rho)-C_L(\rho_{\mathrm{opt}})-
(\mathit{MSE}_m(\rho)-\mathit{MSE}_m(\rho_{\mathrm{opt}}))$ can be bounded by exponential
inequalities given in Lemma 3 of Kneip \cite{Kn94} [the squared norm
$q_\mu(\mathbf{H}_\rho,\mathbf{H}_{\rho_{\mathrm{opt}}})^2$ appearing in
these inequalities can be bounded by $2\mathit{MSE}_m(\rho)$]. These results
lead to
\begin{eqnarray}
\label{gcvproof1}
C_L(\rho)-C_L(\rho_{\mathrm{opt}})&=&\mathit{MSE}_m(\rho)-\mathit{MSE}_m(\rho_{\mathrm{opt}})\nonumber\\[-8pt]\\[-8pt]
&&{}+\eta_{\rho;m}^{[1]}
n^{-{1}/{2}}\mathit{MSE}_m(\rho)^{{1}/{2}},\nonumber
\\
\label{gcvproof2}
\mathit{ASE}_m(\rho)-\mathit{ASE}_m(\rho_{\mathrm{opt}})&=&
\mathit{MSE}_m(\rho)-\mathit{MSE}_m(\rho_{\mathrm{opt}})\nonumber\\[-8pt]\\[-8pt]
&&{}+\eta_{\rho;m}^{[2]}
n^{-{1}/{2}}\mathit{MSE}_m(\rho)^{{1}/{2}},\nonumber
\\
\label{gcvproof3}
\frac{1}{n}\|
\mathbf{Y}-\mathbf{H}_\rho\mathbf{Y}\|^2&=&\sigma^2_\varepsilon+\mathit{MSE}_m(\rho
_{\mathrm{opt}})+\eta_{\rho;m}^{[3]}n^{-{1}/{2}},
\end{eqnarray}
where $\eta_{\rho;m}^{[s]}$ are random variables satisfying
$\sup_{\rho>0}|\eta_{\rho;m}^{[s]}|=O_P(1)$, $s=1,2,3$. By our
assumptions and the arguments used in the proof of Theorem \ref{thm1} we can
infer that $n^{-1}\operatorname{Tr}(\mathbf{H}_\rho)=
O_P([n\rho^{{1}/({2m+2q+1})}]^{-1})= o_P(1)$ for all
$\rho\in[n^{-2m+\delta},\infty)$ as $n\rightarrow\infty$.
Furthermore, there exists a constant $D<\infty$ such that
$n^{-1}\operatorname{Tr}(\mathbf{H}_\rho)\le D\cdot \mathit{MSE}_m(\rho)
=O_P(\rho+[n\rho^{{1}/({2m+2q+1})}]^{-1})$. Together with
(\ref{gcvproof3}) a Taylor expansion of $\mathit{GCV}_m(\rho)$ with respect
to $n^{-1}\operatorname{Tr}(\mathbf{H}_\rho)$ then yields
\begin{eqnarray}\label{gcvproof4}
\mathit{GCV}_m(\rho)&=&\frac{1}{n}\|
\mathbf{Y}-\mathbf{H}_\rho\mathbf{Y}\|^2+2\frac{1}{n}\|
\mathbf{Y}-\mathbf{H}_\rho\mathbf{Y}\|^2\frac{\operatorname{Tr}(\mathbf{H}_\rho)}{n}\nonumber\\
&&{} +
\eta_{\rho;m}^{[4]}\biggl(\frac{\operatorname{Tr}(\mathbf{H}_\rho)}{n}\biggr)^2\\
&=&C_L(\rho)+\eta_{\rho;m}^{[5]}\bigl(n^{-{1}{2}}+\mathit{MSE}_m(\rho)\bigr)\frac
{\operatorname{Tr}(\mathbf{H}_\rho)}{n},\nonumber
\end{eqnarray}
where again $\eta_{\rho;m}^{[s]}$ are random variables with
$\sup_{\rho>n^{-2m+\delta}}|\eta_{\rho;m}^{[s]}|=O_P(1)$, $s=4,5$.
Together with $\mathit{MSE}_m(\rho_{\mathrm{opt}})=O_P(n^{-{2m+2q+1}/({2m+2q+2})})$,
Relation (\ref{gcvrho}) now is an immediate consequence of
(\ref{gcvproof1})--(\ref{gcvproof4}).

Since Lemma 3 of Kneip \cite{Kn94} provides exponential inequalities,
it is easily verified that
uniform bounds similar to (\ref{gcvproof1})--(\ref{gcvproof4}) hold
for all
$\rho\in[n^{-2m+\delta},\infty)$ and \textbf{all} $m=1,\ldots,M_n$, if
$\eta_{\rho;m}^{[s]}$ are
replaced by $\tilde\eta_{\rho;m}^{[s]}\cdot\log M_n$,
$s=1,\ldots,5$. Then $\sup_{\rho>n^{-2m+\delta},m=1,\ldots,M_n}|\tilde
\eta_{\rho;m}^{[s]}|=O_P(1)$, $s=1,\ldots,5$.
The proof of (\ref{gcvrhom}) then follows the arguments used above.

\subsection[Proof of Theorem 4]{Proof of Theorem \protect\ref{thm4}}

Consider the following decomposition:
\[
\widehat{\bolds{\alpha}}_{\mathbf{W}} - \widehat{\bolds
{\alpha}} = \biggl( \frac{1}{n p^{2}} \mathbf{X}^{\tau} \mathbf{X} +
\frac{\rho}{p} \mathbf{A}_{m} \biggr)^{-1}\frac{1}{n p} \bolds
{\delta}^{\tau} \mathbf{Y} +\mathbf{S}\biggl[ \frac{1}{n p} \mathbf
{W}^{\tau} \mathbf{Y} \biggr],
\]
where
\begin{eqnarray*}
\mathbf{S} &:=& \biggl( \frac{1}{n p^{2}} \mathbf{X}^{\tau} \mathbf{X} +
\frac{\rho}{p}
\mathbf{A}_{m} + \mathbf{T} \biggr)^{-1} - \biggl( \frac{1}{n p^{2}}
\mathbf{X}^{\tau} \mathbf{X} + \frac{\rho}{p} \mathbf{A}_{m}
\biggr)^{-1},\qquad
\\
\mathbf{T} &:=& \mathbf{R} - \frac{\widehat{\sigma}_{\delta}^{2} - \sigma
_{\delta}^{2}}{p^{2}} \mathbf{I}_{p}
\end{eqnarray*}
and where $\bolds{\delta}$ is the $n\times p$ matrix
with generic element $\delta_{ij}-\overline{\delta}_j$, $i=1,\ldots,n$,
$j=1,\ldots,p$ and the matrix $\mathbf{R}$ is defined in (\ref{debruitage}).
Thus one obtains
\begin{eqnarray}\label{diffalpha}
\| \widehat{\bolds{\alpha}}_{\mathbf{W}} - \widehat
{\bolds{\alpha}} \|_{\Gamma_{n,p}}
&\leq&\biggl\| \biggl( \frac{1}{n p^{2}} \mathbf{X}^{\tau} \mathbf{X} +
\frac{\rho}{p} \mathbf{A}_{m} \biggr)^{-1} \frac{1}{n p} \bolds
{\delta}^{\tau}
\mathbf{Y} \biggr\|_{\Gamma_{n,p}}\nonumber\\[-8pt]\\[-8pt]
&&{} + \biggl\| \mathbf{S} \biggl(\frac
{1}{n p} \mathbf{W}^{\tau}
\mathbf{Y} \biggr)\biggr\|_{\Gamma_{n,p}}.
\nonumber
\end{eqnarray}
Note that $\mathbb{E}_\varepsilon((\frac{1}{n p^{2}}
\mathbf{X}^{\tau} \mathbf{X} + \frac{\rho}{p} \mathbf{A}_{m}
)^{-1}\frac{1}{n p} \bolds{\delta}^{\tau} \mathbf{Y})=0$ ,
whereas with assumptions (A.1)~and~(A.2)
\begin{eqnarray*}
&&\mathbb{E}_\varepsilon\biggl(\biggl\| \biggl( \frac{1}{n p^{2}} \mathbf
{X}^{\tau} \mathbf{X} + \frac{\rho}{p} \mathbf{A}_{m} \biggr)^{-1}\frac
{1}{n p} \bolds{\delta}^{\tau} \mathbf{Y} \biggr\|^2_{\Gamma
_{n,p}}\biggr)\\
&&\qquad=\mathbb{E}_\varepsilon\biggl(\frac{1}{n^2p}\mathbf{Y}^{\tau}\bolds
{\delta}\biggl( \frac{1}{n p} \mathbf{X}^{\tau} \mathbf{X} + \rho\mathbf
{A}_{m} \biggr)^{-1}\frac{1}{n p} \mathbf{X}^{\tau} \mathbf{X}\biggl(
\frac{1}{n p} \mathbf{X}^{\tau} \mathbf{X} + \rho\mathbf{A}_{m}
\biggr)^{-1}\bolds{\delta}^{\tau}\mathbf{Y}\biggr)\\
&&\qquad= O_P\biggl(\frac{\sigma_\delta^2}{n p}\operatorname{Tr}\biggl(\biggl( \frac{1}{n p}
\mathbf{X}^{\tau} \mathbf{X} + \rho\mathbf{A}_{m}
\biggr)^{-1}\biggr)\biggr).
\end{eqnarray*}
This leads with the properties of the eigenvalues of $(
\frac{1}{n p} \mathbf{X}^{\tau} \mathbf{X}
+ \rho\mathbf{A}_{m} )^{-1}$ to
%
\begin{equation}\label{step1}
\biggl\| \biggl( \frac{1}{n p^{2}} \mathbf{X}^{\tau} \mathbf{X} +
\frac{\rho}{p} \mathbf{A}_{m} \biggr)^{-1} \frac{1}{n p} \bolds
{\delta}^{\tau} \mathbf{Y} \biggr\|_{\Gamma_{n,p}} = O_P\biggl(\frac{1}{(np\rho)^{1/2}}\biggr).
\end{equation}
The next step consists in studying the behavior of the matrix
$\mathbf{R}$ defined in (\ref{debruitage}). Its generic term is
$R_{r,s}=\frac{1}{n p^2}\sum_{i=1}^n(X_i(t_r)-\overline{X}(t_r))(\delta
_{is}-\overline{\delta}_s)+(X_i(t_s)-\overline{X}(t_s))
(\delta_{ir}-\overline{\delta}_r)+(\delta_{ir}-\overline{\delta
}_r)(\delta_{is}-\overline{\delta}_s)-\sigma_\delta^2I{[r=s]}$,
for $r,s=1,\ldots,p$, so that for any $\mathbf{u}\in\mathbb{R}^p$ such
that $\|\mathbf{u}\|=1$ one has
$\|\mathbb{E}_\varepsilon(\mathbf{R}\mathbf{u})\|=O_P(\frac
{1}{np^2})$ whereas it is easy to see that with assumptions (A.1) and
(A.2) and (\ref{conddelta}), $\mathbb{E}_\varepsilon(\|\mathbf{R}\mathbf
{u}\|^2)=O_P(\frac{1}{np^2})$ and then
$\|\mathbf{R}\|=O_P(\frac{1}{n^{1/2}p})$. Now to derive an
upper bound for the norm of the matrix $\mathbf{T}$,
we use the convergence result given in Gasser, Sroka and Jennen-Steinmetz \cite
{GaSrSt86} which in our framework implies that
$
\widehat{\sigma}_{\delta}^{2} = \sigma_{\delta}^{2} + O_{P} ( \frac
{1}{n^{1/2} p} ).
$
Together with the order of $\| \mathbf{R} \|$ this yields
%
\begin{equation}\label{step2}
\| \mathbf{T} \| = O_{P} \biggl( \frac{1}{n^{1/2} p} \biggr).
\end{equation}
For the second term in (\ref{diffalpha}) we consider at first
its Frobenius norm. We have
\begin{eqnarray*}
&&\biggl\| \mathbf{S} \biggl(\frac{1}{n p} \mathbf{W}^{\tau} \mathbf{Y}
\biggr)\biggr\|_{F} \\
&&\qquad\leq\frac{1}{p^{1/2}}\biggl\| \biggl[\biggl( \frac{1}{n p^2} \mathbf
{X}^{\tau} \mathbf{X} + \frac{\rho}{p} \mathbf{A}_{m} + \mathbf{T}
\biggr)^{-1}- \biggl( \frac{1}{n p^2} \mathbf{X}^{\tau} \mathbf{X} + \frac{\rho
}{p} \mathbf{A}_{m} \biggr)^{-1}\biggr] \\
&&\hspace*{194pt}{} \times\biggl(\frac{1}{n^2 p^2} \mathbf{W}^{\tau} \mathbf
{Y}\mathbf{Y}^\tau\mathbf{W}\biggr)^{1/2}\biggr\|_F \\
&&\qquad\leq\frac{1}{p^{1/2}}\biggl\| \biggl( \frac{1}{n p^2} \mathbf{X}^{\tau}
\mathbf{X} + \frac{\rho}{p} \mathbf{A}_{m} \biggr)^{-1}\frac{1}{n p}
\mathbf{W}^{\tau} \mathbf{Y}\biggr\|^2 \|\mathbf{T}\| \biggl\|
\frac{1}{n p} \mathbf{W}^{\tau} \mathbf{Y}\biggr\|^{-1},
\end{eqnarray*}
where the second inequality comes from the first inequality
in Demmel \cite{De92}. Note that with assumptions (A.2) and (A.5), for
every $\delta>0$, there is a positive constant such that $p^{1/2}\|
\mathbb{E}_\varepsilon( \frac{1}{n p} \mathbf{W}^{\tau} \mathbf{Y}
)\|$ is greater than this constant with a probability
larger than or equal to $1-\delta$. We also have $\mathbb{E}_\varepsilon
(\| \frac{1}{n p} \mathbf{W}^{\tau} \mathbf{Y} - \mathbb
{E}_\varepsilon( \frac{1}{n p} \mathbf{W}^{\tau} \mathbf{Y}
)\|^2)$, which is of order $\frac{1}{n p}$.
This gives finally when combining (\ref{alphahboundnew}), (\ref{step2})
and the condition on $p$ and $\rho$ as well as assumption (A.2)
%
\begin{equation}\label{thm4step6}
\biggl\| \mathbf{S} \frac{1}{n p} \mathbf{W}^{\tau} \mathbf{Y}\biggr\|
_{\Gamma_{n,p}}^2 = O_P\biggl(\biggl\| \mathbf{S} \biggl(\frac{1}{n p}
\mathbf{W}^{\tau} \mathbf{Y} \biggr)\biggr\|_{F}\biggr) = O_P
\biggl(\frac{1}{n}\biggr),
\end{equation}
which concludes Theorem \ref{thm4} with (\ref{diffalpha}) and
(\ref{step1}).

\subsection[Proof of Theorem 5]{Proof of Theorem \protect\ref{thm5}}

We first prove (\ref{resthm5}). Obviously,
\[
\| \widehat{\alpha}_{\mathbf{W}} - \widehat{\alpha} \|
_{\Gamma_n}^{2} \leq\frac{2}{n} \sum_{i=1}^{n} (\widehat
{d}_{i,\mathbf{W}}-\overline{\widehat{d}}_{\mathbf{W}}
)^{2} + 2\| \widehat{\bolds{\alpha}}_{\mathbf{W}} -
\widehat{\bolds{\alpha}} \|_{\Gamma_{n,p}}^{2},
\]
where
\[
\widehat{d}_{i,\mathbf{W}} = \int_{I} \bigl(\widehat{\alpha
}_{\mathbf{W}} (t) - \widehat{\alpha}(t)\bigr) X_{i} (t) \,dt - \frac
{1}{p} \sum_{j=1}^{p} \bigl(\widehat{\alpha}_{\mathbf{W}} (t_{j}) -
\widehat{\alpha}(t_j)\bigr)X_{i} (t_{j}).
\]
Then, assertion (\ref{resthm4}) implies that (\ref{resthm5}) is a
consequence of
%
\begin{eqnarray}\label{endthm5}
\frac{1}{n} \sum_{i=1}^{n} (\widehat{d}_{i,\mathbf
{W}}-\overline{\widehat{d}}_{\mathbf{W}})^{2} = O_P\biggl(
\frac{1}{n p \rho} + \frac{1}{n}\biggr).
\end{eqnarray}
The proof of (\ref{endthm5}) follows the same structure as
the proof of (\ref{endthm2}). Indeed, we have
\begin{eqnarray}\label{thm5step1}\hspace*{28pt} 
&&\frac{1}{n}\sum_{i=1}^n(\widehat{d}_{i,\mathbf
{W}}-\overline{\widehat{d}}_{\mathbf{W}})^2 
\nonumber\\
&&\qquad\le2 x_{\max}^2 \Biggl(
\sum_{j=1}^p \frac{1}{p} \biggl[\int_{t_j-1/(2p)}^{t_j+1/(2p)} |P'(t)|
+|P'_{\mathbf{W}}| + |r(t)| + |r_{\mathbf{W}}| \,dt\biggr]^2\Biggr)
\nonumber\\[-8pt]\\[-8pt]
&&\qquad\quad{} +2 \frac{1}{p}
\|\widehat{\bolds{\alpha}}_{\mathbf{W}} - \widehat
{\bolds{\alpha}}\|^2\nonumber\\
&&\qquad\quad\hspace*{10.6pt}{}\times \frac{1}{n}\sum_{i=1}^n
\sum_{j=1}^p \int_{t_j-1/(2p)}^{t_j+1/(2p)}\bigl(\bigl(X_i(t)-\overline
{X}(t)\bigr)-\bigl(X_i(t_j)-\overline{X}(t_j)\bigr)\bigr)^2\,dt, \nonumber
\end{eqnarray}
where $P_{\mathbf{W}}(t) =
\sum_{l=0}^{m-1}\frac{t^l}{l!}\widehat{\alpha}_{\mathbf{W}}(0)$,
$r_{\mathbf{W}}(t)
= \int_0^t \frac{(t-u)^{m-1}}{(m-1)!}\widehat{\alpha}_{\mathbf
{W}}(u)\,du$ and $P(t)$ and $r(t)$
are similarly defined for $\widehat{\alpha}$ (see the proof of Theorem
\ref{thm2}).

Replacing the semi-norm $\Gamma_{n,p}$ by the euclidean norm
in (\ref{resthm4}) following the same lines
as the proof of Theorem \ref{thm4}, one can show that
%
\begin{equation}\label{thm5step2}
\frac{1}{p} \|\widehat{\bolds{\alpha}}_{\mathbf{W}} -
\widehat{\bolds{\alpha}} \|^2 = \frac{1}{p}(\widehat{\bolds{\alpha}}_{\mathbf{W}} -
\widehat{\bolds{\alpha}})^\tau(\widehat{\bolds
{\alpha}}_{\mathbf{W}}
- \widehat{\bolds{\alpha}}) = O_P\biggl(\frac{1}{n p \rho
^2} + \frac{1}{n}\biggr),
\end{equation}
which together with assumption (A.2) implies that the second
term on the right-hand side of (\ref{thm5step1}) can
be bounded by $O_P(\frac{p^{-2\kappa}}{n p \rho^2} + \frac
{p^{-2\kappa}}{n})$.


Now the remainder of the proof consists in
studying $\int_0^1\widehat{\alpha}_{\mathbf{W}}^{(m)}(t)^2\,dt$.
Recalling the definition of
$\widehat{\bolds{\alpha}}_{\mathbf{W}}$, we have
\begin{eqnarray*}
&& \frac{1}{n} \biggl\| \mathbf{Y} - \frac{1}{p} \mathbf{W} \widehat
{\bolds{\alpha}}_{\mathbf{W}} \biggr\|^{2} + \frac{\rho}{p}
\widehat{\bolds{\alpha}}_{\mathbf{W}}^\tau\mathbf{P}_m
\widehat{\bolds{\alpha}}_{\mathbf{W}} + \rho\int_{I} \widehat
{\alpha}_{\mathbf{W}}^{(m)} (t)^{2} \,dt - \frac{\widehat{\sigma
}_\delta}{p^2} \widehat{\bolds{\alpha}}_{\mathbf{W}}^\tau
\widehat{\bolds{\alpha}}_{\mathbf{W}}\\
&&\qquad \leq \frac{1}{n} \biggl\| \mathbf{Y} - \frac{1}{p} \mathbf{W}
\widehat{\bolds{\alpha}} \biggr\|^{2} + \frac{\rho}{p} \widehat
{\bolds{\alpha}}^\tau\mathbf{P}_m \widehat{\bolds{\alpha
}} + \rho\int_{I} \widehat{\alpha}^{(m)} (t)^{2} \,dt - \frac{\widehat
{\sigma}_\delta}{p^2} \widehat{\bolds{\alpha}}^\tau \widehat
{\bolds{\alpha}}
\end{eqnarray*}
and then
\begin{eqnarray}\label{thm5ineg1}
&&
\rho\int_{I} \widehat{\alpha}_{\mathbf{W}}^{(m)} (t)^{2} \,dt\nonumber\\
&&\qquad\leq\frac{1}{n} \biggl\| \frac{1}{p} \mathbf{W}(\widehat
{\bolds{\alpha}}_{\mathbf{W}} - \widehat{\bolds{\alpha
}}) \biggr\|^{2} + \frac{2}{n} \biggl\langle\mathbf{Y} - \frac
{1}{p} \mathbf{W} \widehat{\bolds{\alpha}} , \frac{1}{p} \mathbf
{W} \widehat{\bolds{\alpha}} - \frac{1}{p} \mathbf{W} \widehat
{\bolds{\alpha}}_{\mathbf{W}}\biggr\rangle\nonumber\\[-8pt]\\[-8pt]
&&\qquad\quad{} - \frac{\rho}{p} \widehat{\bolds{\alpha}}_{\mathbf{W}}^\tau
\mathbf{P}_m \widehat{\bolds{\alpha}}_{\mathbf{W}} + \frac
{\rho}{p} \widehat{\bolds{\alpha}}^\tau\mathbf{P}_m \widehat
{\bolds{\alpha}}
\nonumber\\
&&\qquad\quad{} + \frac{\widehat{\sigma}_\delta}{p^2} \widehat{\bolds{\alpha
}}_{\mathbf{W}}^\tau \widehat{\bolds{\alpha}}_{\mathbf
{W}} - \frac{\widehat{\sigma}_\delta}{p^2} \widehat{\bolds{\alpha
}}^\tau \widehat{\bolds{\alpha}} + \rho\int_{I} \widehat{\alpha
}^{(m)} (t)^{2} \,dt.\nonumber
\end{eqnarray}
First consider the term $\frac{1}{n} \| \frac{1}{p}
\mathbf{W}(\widehat{\bolds{\alpha}}_{\mathbf{W}} -
\widehat{\bolds{\alpha}}) \|^{2}$.
By (\ref{resthm4}) and (\ref{thm5step2}) we obtain
\begin{eqnarray}\label{step3}
\frac{1}{n} \biggl\| \frac{1}{p} \mathbf{W}(\widehat
{\bolds{\alpha}}_{\mathbf{W}} -
\widehat{\bolds{\alpha}}) \biggr\|^{2}
= O_P\biggl(\frac{1}{np\rho} + \frac{1}{n}\biggr).
\end{eqnarray}
We focus now on the second term in the right-hand side of (\ref
{thm5ineg1}), for which we have the following decomposition:
\begin{eqnarray*}
&& \frac{1}{n} \biggl\langle\mathbf{Y} - \frac{1}{p} \mathbf{W} \widehat
{\bolds{\alpha}} , \frac{1}{p} \mathbf{W} \widehat{\bolds
{\alpha}} - \frac{1}{p} \mathbf{W} \widehat{\bolds{\alpha
}}_{\mathbf{W}}\biggr\rangle\\
&&\qquad = \frac{1}{n} \biggl\langle\frac{1}{p}
\mathbf{X} \bolds{\alpha} - \frac{1}{p} \mathbf{W} \widehat
{\bolds{\alpha}} , \frac{1}{p} \mathbf{W} \widehat{\bolds
{\alpha}} - \frac{1}{p} \mathbf{W} \widehat{\bolds{\alpha
}}_{\mathbf{W}}\biggr\rangle\\
&&\qquad\quad{}  + \frac{1}{n}\biggl\langle \mathbf{d} , \frac{1}{p}
\mathbf{W} \widehat{\bolds{\alpha}} - \frac{1}{p} \mathbf{W}
\widehat{\bolds{\alpha}}_{\mathbf{W}}\biggr\rangle + \frac
{1}{n}\biggl\langle \bolds{\varepsilon} , \frac{1}{p} \mathbf{W}
\widehat{\bolds{\alpha}} - \frac{1}{p} \mathbf{W} \widehat
{\bolds{\alpha}}_{\mathbf{W}}\biggr\rangle.
\end{eqnarray*}
We have
\begin{eqnarray*}
&&\frac{1}{n^{1/2}} \biggl\| \frac{1}{p} \mathbf{X} \bolds{\alpha} -
\frac{1}{p} \mathbf{W} \widehat{\bolds{\alpha}}\biggr\| \\
&&\qquad\leq\frac{1}{n^{1/2}} \biggl\| \frac{1}{p} \mathbf{X} \bolds
{\alpha} - \frac{1}{p} \mathbf{W} \widehat{\bolds{\alpha}} -
\mathbb{E}_\varepsilon\biggl( \frac{1}{p} \mathbf{X} \bolds{\alpha} -
\frac{1}{p} \mathbf{W} \widehat{\bolds{\alpha}} \biggr)\biggr\|\\
&&\qquad\quad{} +
\frac{1}{n^{1/2}} \biggl\| \mathbb{E}_\varepsilon\biggl(\frac{1}{p} \mathbf
{X} \bolds{\alpha} - \frac{1}{p} \mathbf{W} \widehat{\bolds
{\alpha}}\biggr)\biggr\|.
\end{eqnarray*}
Some straightforward calculations and previous results lead to
$\frac{1}{n^{1/2}} \| \frac{1}{p} \mathbf{X} \bolds{\alpha
} - \frac{1}{p} \mathbf{W} \widehat{\bolds{\alpha}}
- \mathbb{E}_\varepsilon( \frac{1}{p} \mathbf{X} \bolds{\alpha}
- \frac{1}{p} \mathbf{W} \widehat{\bolds{\alpha}} )\|
=O_P( (1/n \rho^{1/(2m+2q+1)})^{1/2} +1/p^{1/2})$ whereas
$\|\mathbb{E}_\varepsilon( \frac{1}{p} \mathbf{X} \bolds
{\alpha} - \frac{1}{p} \mathbf{W} \widehat{\bolds{\alpha}}
)\| = O_P(\rho^{1/2} +p^{-\kappa})$.
This finally leads with the Cauchy--Schwarz inequality to
\begin{eqnarray}\label{thm5step5}
&&  \frac{1}{n} \biggl\langle\frac{1}{p} \mathbf{X}
\bolds{\alpha} - \frac{1}{p} \mathbf{W} \widehat{\bolds{\alpha}}
, \frac{1}{p} \mathbf{W} \widehat{\bolds{\alpha}} - \frac{1}{p}
\mathbf{W} \widehat{\bolds{\alpha}}_{\mathbf{W}}\biggr\rangle
\nonumber\\
&&\qquad = O_P\biggl(\biggl(\biggl(\frac{1}{n \rho^{1/2m+2q+1}}\biggr)^{1/2} +
\frac{1}{p^{1/2}} + \rho^{1/2} + p^{-\kappa}\biggr)\\
&&\hspace*{151.5pt}{}\times \biggl(\frac{1}{(n p
\rho)^{1/2}} + \frac{1}{n^{1/2}}\biggr)\biggr).\nonumber
\end{eqnarray}
Using again the Cauchy--Schwarz inequality and (\ref{step3})
we have
%
\begin{equation}\label{thm5step6}
\frac{1}{n}\biggl\langle \mathbf{d} , \frac{1}{p} \mathbf{W} \widehat
{\bolds{\alpha}} -
\frac{1}{p} \mathbf{W} \widehat{\bolds{\alpha}}_{\mathbf
{W}}\biggr\rangle= O_P\biggl( \frac{p^{-\kappa}}{(np\rho)^{1/2}}
+ \frac{p^{-\kappa}}{n^{1/2}}\biggr).
\end{equation}
The last term is such that
\begin{eqnarray*}
&&\frac{1}{n}\bolds{\varepsilon}^\tau\biggl(\frac{1}{p} \mathbf{W}
(\widehat{\bolds{\alpha}} - \widehat{\bolds{\alpha
}}_{\mathbf{W}})\biggr) \\
&&\qquad = \frac{1}{n}\bolds{\varepsilon}^\tau\biggl(\frac{1}{p} \mathbf{W}
\biggl(\frac{1}{n p^2}\mathbf{X}^\tau\mathbf{X} + \frac{\rho}{p} \mathbf
{A}_m\biggr)^{-1} \bolds{\delta}^\tau\mathbf{Y} \biggr) + \frac
{1}{n}\bolds{\varepsilon}^\tau\biggl(\frac{1}{p} \mathbf{W} \mathbf
{S}\biggl(\frac{1}{n p } \mathbf{W}^\tau\mathbf{Y}\biggr)\biggr).
\end{eqnarray*}
Using the same developments as above and using assumptions
(A.1) and~(A.2)
we obtain that $\frac{1}{n}\bolds{\varepsilon}^\tau(\frac{1}{p}
\mathbf{W}
(\frac{1}{n p^2}\mathbf{X}^\tau\mathbf{X} + \frac{\rho}{p} \mathbf
{A}_m)^{-1} \bolds{\delta}^\tau\mathbf{Y}
)=O_P(\frac{1}{n p^{1/2}\rho^{1/2}})$ while\break $\frac
{1}{n}\bolds{\varepsilon}^\tau(\frac{1}{p} \mathbf{W}
\mathbf{S}\times(\frac{1}{n p } \mathbf{W}^\tau\mathbf{Y}))
= O_P(\frac{1}{n})$. This finally leads to
%
\begin{equation}\label{thm5step7}
\frac{1}{n}\bolds{\varepsilon}^\tau\biggl(\frac{1}{p} \mathbf{W}
(\widehat{\bolds{\alpha}} - \widehat{\bolds{\alpha
}}_{\mathbf{W}})\biggr) = O_P\biggl( \frac{1}{n p^{1/2}\rho
^{1/2}} + \frac{1}{n}\biggr).
\end{equation}
Finally using the same arguments as in the proof of Theorem
\ref{thm2},
assertion (\ref{endthm5}) is a consequence of (\ref{thm5step1}), (\ref
{alphahboundnew})
and (\ref{intalbound}) as well as the bounds obtained in (\ref
{thm5step2})--(\ref{thm5step7}) and the conditions on $n$, $p$ and $\rho$.

It remains to show (\ref{resthm6}).
The proof follows the same lines as the proof of Theorem \ref{thm3}. We
have the following relation:
\begin{eqnarray*}
&& \| \widehat{\alpha}_{\mathbf{W}} - \widehat{\alpha} \|
_{\Gamma_{n}}^{2}\\
&&\qquad =
\| \widehat{\alpha}_{\mathbf{W}} - \widehat{\alpha} \|
_{\Gamma}^{2} + \sum_{r=1}^{\infty} \sum_{s=1}^{\infty}
\widetilde\alpha_{\mathbf{W},r} \widetilde\alpha_{\mathbf{W},s} \Biggl(
\frac{1}{n}\sum_{i=1}^n\tau_{ri}\tau_{si}-\lambda_rI(r=s)\Biggr)+O_P(n^{-1}),
\end{eqnarray*}
with $\widetilde\alpha_{\mathbf{W},r} = \langle\zeta_r ,
\widehat{\alpha}_{\mathbf{W}} - \widehat{\alpha} \rangle$.
Using the Cauchy--Schwarz inequality as in (\ref{gammaalphacauchy}),
the remainder of the proof
consists in showing that $\|\widehat{\alpha}_{\mathbf{W}}-\widehat
{\alpha}\|=O_P(1)$. This
is obtained by using the bounds obtained in the proof of (\ref
{resthm5}) and following the same lines of argument as
for showing (\ref{alphahboundnew}).

\printaddresses


\begin{thebibliography}{99}

\bibitem{AnCaEsVi04}
\textsc{Aneiros-Perez, G., Cardot, H., Estevez-Perez, G.}
and \textsc{Vieu, P.} (2004). Maximum ozone concentration forecasting by
functional nonparametric approaches. \textit{Environmetrics}
\textbf{15} 675--685.

\bibitem{Bo00}
\textsc{Bosq, D.} (2000). \textit{Linear Processes in Function
Spaces. Lecture Notes in Statist.} \textbf{149}. Springer, New York.
\MR{1783138}

\bibitem{Ca00}
\textsc{Cardot, H.} (2000).
Nonparametric estimation of smoothed principal components analysis
of sampled noisy functions.
\textit{J. Nonparametr. Statist.} \textbf{12} 503--538.
\MR{1785396}

\bibitem{CaHa06}
\textsc{Cai, T. T.} and \textsc{Hall, P.} (2006). Prediction in
functional linear regression.
\textit{Ann. Statist.} \textbf{34} 2159--2179.
\MR{2291496}

\bibitem{CaCrKnSa07}
\textsc{Cardot, H., Crambes, C., Kneip, A.} and \textsc{Sarda, P.} (2007).
Smoothing splines estimators in functional linear regression with
errors-in-variables.
\textit{Comput. Statist. Data Anal.} \textbf{51} 4832--4848.
\MR{2364543}

\bibitem{CaCrSa06} \textsc{Cardot, H., Crambes, C.} and \textsc{Sarda, P.}
(2007). Ozone
pollution forecasting. In \textit{Statistical Methods for
Biostatistics and Related Fields} (W. H\"ardle, Y. Mori and P. Vieu,
eds.) 221--244. Springer, New York.
\MR{2376412}

\bibitem{CaFeSa03}
\textsc{Cardot, H., Ferraty, F.} and \textsc{Sarda, P.} (2003). Spline
estimators for the functional linear model. \textit{Statist. Sinica}
\textbf{13} 571--591.
\MR{1997162}

\bibitem{CaMaSa07}
\textsc{Cardot, H., Mas, A.} and \textsc{Sarda, P.} (2007). CLT in
functional linear regression models.
\textit{Probab. Theory Related Fields} \textbf{138} 325--361.
\MR{2299711}

\bibitem{ChMuWa03}
\textsc{Chiou, J. M., M\"uller, H. G.} and \textsc{Wang, J. L.} (2003).
Functional quasi-likelihood regression models with smoothed
random effects.
\textit{J. Roy. Statist. Soc. Ser. B} \textbf{65} 405--423.
\MR{1983755}

\bibitem{CuFeFr02}
\textsc{Cuevas, A., Febrero, M.} and \textsc{Fraiman, R.} (2002).
Linear functional regression: The case of a fixed design and functional
response.
\textit{Canadian J. Statistics} \textbf{30} 285--300.
\MR{1926066}


\bibitem{De92}
\textsc{Demmel, J.} (1992). The componentwise distance to the
nearest singular matrix. \textit{SIAM J. Matrix Anal. Appl.}
\textbf{13} 10--19.
\MR{1146648}

\bibitem{EiMa96}
\textsc{Eilers, P. H.} and \textsc{Marx, B. D.}
(1996). Flexible smoothing
with B-splines and penalties. \textit{Statist. Sci.}
\textbf{11} 89--102.
\MR{1435485}

\bibitem{Eu88}
\textsc{Eubank, R. L.} (1988). \textit{Spline Smoothing and
Nonparametric Regression}. Dekker, New York.
\MR{0934016}

\bibitem{FeVi06}
\textsc{Ferraty, F.} and \textsc{Vieu, P.} (2006). \textit{Nonparametric
Functional Data Analysis: Methods, Theory, Applications and
Implementations}. Springer, London.
\MR{2229687}

\bibitem{Fu87}
\textsc{Fuller, W. A.} (1987). \textit{Measurement Error Models}.
Wiley, New York.
\MR{0898653}


\bibitem{GaSrSt86}
\textsc{Gasser, T., Sroka, L.} and \textsc{Jennen-Steinmetz, C.}
(1986). Residual variance and residual pattern in nonlinear regression.
\textit{Biometrika} \textbf{3} 625--633.
\MR{0897854}

\bibitem{GoVanLo80}
\textsc{Golub, G. H.} and \textsc{Van Loan, C. F.} (1980). An analysis
of the total least squares problem.
\textit{SIAM J. Numer. Anal.}
\textbf{17} 883--893.
\MR{0595451}

\bibitem{HaHo07}
\textsc{Hall, P.} and \textsc{Horowitz, J. L.} (2007). Methodology and
convergence rates for functional linear regression.
\textit{Ann. Statist.} To appear.
\MR{2332269}

\bibitem{HeMuWa00}
\textsc{He, G., M\"uller, H.-G.} and \textsc{Wang, J. L.} (2000).
Extending correlation and regression from multivariate to functional
data. In \textit{Asymptotics in Statistics and Probability} (M. L. Puri,
ed.) 301--315. VSP, Leiden. 

\bibitem{Kn94}
\textsc{Kneip, A.} (1994). Ordered linear smoothers.
\textit{Ann. Statist.} \textbf{22} 835--866.
\MR{1292543}

\bibitem{Li06}
\textsc{Li, Y.} and \textsc{Hsing, T.} (2006). On rates of convergence in
functional linear regression.
\textit{J. Mulitivariate Anal.}
Published online DOI:
\href{http://dx.doi.org/10.1016/j.jmva.2006.10.004}{10.1016/j.jmva.2006.10.004}.
\MR{2392433}

\bibitem{MaEi99}
\textsc{Marx, B. D.} and \textsc{Eilers, P. H.} (1999).
Generalized linear regression on sampled signals and curves:
A $P$-spline approach. \textit{Technometrics} \textbf{41} 1--13.

\bibitem{MuSt05}
\textsc{M\"uller, H.-G.} and \textsc{Stadtm\"uller, U.} (2005).
Generalized functional linear models. \textit{Annn. Statist.}
\textbf{33} 774--805.
\MR{2163159}

\bibitem{RaDa91}
\textsc{Ramsay, J. O.} and \textsc{Dalzell, C. J.} (1991). Some tools for
functional data analysis. \textit{J. Roy. Statist.
Soc. Ser. B} \textbf{53} 539--572.
\MR{1125714}

\bibitem{RaSi02}
\textsc{Ramsay, J. O.} and \textsc{Silverman, B. W.} (2002).
\textit{Applied Functional Data Analysis}. Springer, New York.
\MR{1910407}

\bibitem{RaSi05}
\textsc{Ramsay, J. O.} and \textsc{Silverman, B. W.} (2005).
\textit{Applied Functional Data Analysis}, 2nd ed. Springer, New York.
\MR{2168993}

\bibitem{St82}
\textsc{Stone, C. J.} (1982). Optimal global rates of convergence
for nonparametric regression. \textit{Ann. Statist.}
\textbf{10} 1040--1053.
\MR{0673642}

\bibitem{Ut83}
\textsc{Utreras, F.} (1983). Natural spline functions, their
associated eigenvalue problem. \textit{Numer. Math.}
\textbf{42} 107--117.
\MR{0716477}

\bibitem{VanHuVa91}
\textsc{Van Huffel, S.} and \textsc{Vandewalle, J.} (1991).
\textit{The Total Least Squares Problem: Computational Aspects and Analysis}.
SIAM, Philadelphia.
\MR{1118607}

\bibitem{Wa77}
\textsc{Wahba, G.} (1977). Practical approximate solutions to
linear operator equations when the data are noisy. \textit{SIAM J.
Numer. Anal.} \textbf{14} 651--667.
\MR{0471299}

\bibitem{Wa90}
\textsc{Wahba, G.} (1990). \textit{Spline Models for
Observational Data}. SIAM, Philadelphia.
\MR{1045442}\

\bibitem{YaMuWa05}
\textsc{Yao, F., M\"uller, H.-G.} and \textsc{Wang, J. L.} (2005).
Functional data analysis for sparse longitudinal data.
\textit{J. Amer. Statist. Assoc.} \textbf{100} 577--590.
\MR{2160561}

\end{thebibliography}
\end{document}